\documentclass[reqno,9pt]{amsart}
\usepackage{amsmath,amsxtra,amssymb,latexsym, amscd,amsthm}
\usepackage[unicode]{hyperref}
\usepackage{array,tabularx,longtable,multicol,indentfirst,fancybox,color}
\usepackage{graphicx}
\usepackage{multicol}
\usepackage{mathrsfs}

\advance\hoffset-1truecm\relax
\newtheorem{proposition}{Proposition}[section]
\newtheorem{theorem}{Theorem}[section]
\newtheorem{definition}{Definition}
\newtheorem{lemma}{Lemma}[section]
\newtheorem{corollary}{Corollary}[section]
\newtheorem{remark}{Remark}

\numberwithin{equation}{section}

\def\R{\mathbb{R}}

\def\vp{\varphi}

\def\ga{\gamma}

\def\lv{\lVert}

\begin{document}
	\title[Composition structure of polyconvolution associated with Kontorovich-Lebedev transform and Fourier]{Composition structure of polyconvolution associated with\\ index Kontorovich-Lebedev transform and Fourier integrals}
	\date{15-July-2025, Accepted by \textsc{Integral Transforms and Special Functions}}	
	\author[Trinh Tuan]{Trinh Tuan}

\maketitle	
\begin{center}
	Department of Mathematics, Faculty of Nature Sciences, Electric Power University,\\ 235-Hoang Quoc Viet Rd.,  Hanoi, Vietnam.\\
	E-mail: \texttt{tuantrinhpsac@yahoo.com}
\end{center}	
\begin{abstract}
Using Kakichev's classical concept and extending Yakubovich-Britvina's approach (\textit{Results. Math.} 55(1-2):175-197, 2009) and (\textit{Integral Transforms Spec. Funct.} 21(4):259--276, 2010) for setting up Kontorovich-Lebedev convolution operators, this paper proposes a new polyconvolution structure associated with the KL-transform and Fourier integrals. Our main contributions include demonstrating a one-dimensional Watson-type transform, providing necessary and sufficient conditions for this transform to serve as unitary on $L_2(\mathbb{R}_+)$, and inferring its inverse operator in symmetric form. The existence of this structure over specific function spaces and its connection with previously known convolutions are pointed out. Establish the Plancherel-type theorem, prove the convergence in the mean-square sense in $L_2(\mathbb{R}_+)$, and prove the boundedness of dual spaces via Riesz-Thorin's theorem. Derives new weighted $L_p$-norm inequalities and boundedness in a three-parametric family of Lebesgue spaces. These theoretical findings are applied to solve specific classes of the Toeplitz-Hankel equation, providing a priori estimations based on the established conditions for $L_1$ solvability.

\vskip 0.3cm

	\end{abstract}

\section{Introduction}
\subsection{Context of the problem} The Kontorovich–Lebedev transform (abbreviated as KL-transform) is an integral transform that uses a Macdonald function (modified Bessel function of the second kind) with an imaginary index as its kernel \cite{bateman1954}. Unlike other Bessel function transforms, such as the Hankel transform \cite{GlaeskeHJ2006}, this transform involves integrating over the index of the function rather than its argument. KL-transform was introduced for the first time in \cite{KLtransform1938} to solve certain
boundary-value problems of mathematical physics. It arises naturally when the method of separation of variables is used to solve the Laplace equation in cylindrical coordinates for wedge-shaped domains. For more detailed information on the Kontorovich–Lebedev transformation and related issues, we refer the reader
to \cite{GlaeskeHJ2006,Yakubovich1996index,Yakubovich2008progress,yakubovich1994theory,Zayed1996handbook,Zemanian1975kontorovich}.
Let $f, h$ be functions defined on domain $\mathbb{R}_{+} \equiv [0, \infty)$. The following double integral, which we call as the convolution of $f$ and $h$, is denoted by $(f * h)(x)$, and we define it as in \cite{Yakubovich1996index} by
$$
(f * h)(x)=\frac{1}{2 x} \int_0^{\infty} \int_0^{\infty} e^{-\frac{1}{2}(x \frac{u^2+y^2}{u y}+\frac{y u}{x})} f(u) h(y) d u d y,\ \text{with}\ x>0.
$$
This convolution has a close relationship with the KL-transform \cite{Sneddon1972}
$
K_{i \tau}[f]=\int_0^{\infty} K_{i \tau}(x) f(x) d x$ with $ \tau\in \mathbb{R}_{+}
$, which contains, as its kernel, the modified Bessel function of the second kind $K_\nu(z)$, or Macdonald's function \cite{Erdelyi1953higher} with a purely imaginary index $\nu=i \tau$. The function $K_\nu(z)$ satisfies the second order differential equation
$$
z^2 \frac{d^2 u}{d z^2}+z \frac{d u}{d z}-\left(z^2+\nu^2\right) u=0
$$
for which it is the solution that remains bounded as $z$ tends to infinity on the real line. 
One of the outstanding results concerning constructing convolution operators, which is associated with the KL-transform and Fourier transform, is as follows: If denoted by $\kappa(t, \Theta ; \tau)=K_0\left(\sqrt{t^2+\Theta^2+2 t \Theta \cosh \tau}\right)$, based on the reciprocities for the product of modified Bessel functions $K_{i x}$ with
$
\int_0^{\infty} K_{i x}(\Theta) K_{i x}(t) \cos x \tau d x =\frac{\pi}{2} K_0\left(\sqrt{t^2+\Theta^2+2 t \Theta \cosh \tau}\right)$
and $ K_{i x}(\Theta) K_{i x}(t)  =\int_0^{\infty} K_0\left(\sqrt{t^2+\Theta^2+2 t \Theta \cosh \tau}\right) \cos x \tau d \tau.$
Obviously, this kernel is symmetric by variables $t$ and $\Theta$, this means that $\kappa(t, \Theta ; \tau)=\kappa(\Theta, t ; \tau)$, two convolutions 
$
(f * g)_1(t)=\frac{1}{\pi t} \int_0^{\infty} \int_0^{\infty} f(\tau) g(\Theta) \kappa(t, \Theta ; \tau)\ d \tau d \Theta$ and $
(f * g)_2(t)=\int_0^{\infty} \int_0^{\infty} f(\tau) g(\Theta) \kappa(\tau, \Theta ; t)\ d \tau d \Theta,
$
with the same kernel  $\kappa(t, \Theta ; \tau)$, have been introduced and studied by  Yakubovich and Britvina in \cite{Yakubovich2009Britvina,yakubovich2010Britvina}.  
The essence of an integral transform lies in the properties of the associated kernel functions involved. Therefore, the difficulty in constructing convolutions involving the KL-transform depends entirely on the properties and asymptotic behavior of the Macdonald function \cite{glaeske1986MathZ}. When more than two transformations act on the convolution, we get the terminology of generalized convolution or polyconvolution.  
It is well established that there are two principal methodologies for constructing convolution operations for integral transforms. The first approach involves the creation of a generalized shift operator, also referred to as a generalized translation or displacement operator. In this framework, the standard translation operator appearing in the classical convolution is replaced by its generalized counterpart, thereby yielding a generalized convolution structure adapted to the specific integral transform under consideration.
The second approach originates from the work of Valentin A. Kakichev, whose constructing method is based on the factorization identity. This technique facilitates the construction of convolution formulas for a broad class of integral transforms and is applicable in a highly general setting, as initially presented in \cite{kakichev1967convolution}. Afterward, in his lectures \cite{Kakichev1997}, Kakichev generalized this approach and introduced the concept of polyconvolution.

Being directly influenced by the above is also the driving motivation for this work. We have found a new structure of polyconvolution (generalized convolution) related to the KL-transform. Its component structure aligns with Kakichev's classical concept and has certain connections with the previously known convolutions of Sneddon and Yakubovich-Britvina.  For this reason, our setup will be heavily influenced by the approach in \cite{Yakubovich2009Britvina,yakubovich2010Britvina}. Besides, the existence of this polyconvolution is fully validated via theorems in this article, and its component structure will lead to many intriguing properties, such as the Plancherel theory and the Watson transform in 1-dimensional, among others. Our work has enriched the generalized convolution theory of the Kontorovich–Lebedev integral transform and associated problems.
\subsection{Some notations and auxiliary results}
We briefly recall essential definitions and previous results from \cite{Yakubovich1996index,Sneddon1972,Kakichev1997,Prudnikov1986Marichev,Titchmarsh1986}. The following notations and conventions shall be deployed throughout the article
\begin{definition}[Kakichev's concept \cite{Kakichev1997}]\label{dn1}
Let $U_i (X_i)$ be the space of linear functions that can differ on the same field and $\mathcal{V} (Y)$ an algebra. Consider the integral transforms $\mathcal{K}_i: U_i (X_i) \to \mathcal{V}(Y)$ with $i=1,\dots,n+1.$ Then polyconvolution (generalized convolution) of functions $f_1 \in U_1 (X_1), f_2 \in U_2 (X_2),\dots, f_n \in U_n (X_n)$ with weighted $\ga$ for the integral transforms $\mathcal{K}_i$ is a multi-linear operator defined by $\textasteriskcentered: \prod\limits_{i=1}^n U_i (X_i)\to \mathcal{V}(Y),$ denoted by mapping $\underset{\ga}{*} (f_1 , f_2 , \dots, f_n)$ such that the following factorization property is valid
$\mathcal{K}_{n+1} \big[ \underset{\ga}{*} (f_1 , f_2 , \dots, f_n)\big](y)= \ga(y)\prod\limits_{i=1}^n (\mathcal{K}_i f_i)(y).$ 
Here, we always assume that for any integral transformation $\mathcal{K}_i$, there always exists a corresponding inverse transform.
Notice that multiplication on the right-hand side is understood as multiplication in algebra $\mathcal{V}(Y)$. And
$
\underset{\ga}{*} (f_1 , f_2 , \dots, f_n)(x)=\mathcal{K}^{-1}_{n+1} \big[\ga(y)\prod\limits_{i=1}^n (\mathcal{K}_i f_i)(y) \big](x),
$
where $\mathcal{K}^{-1}_{n+1}$ is the inverse transform of $\mathcal{K}_{n+1}$.
\end{definition}
The Fourier transforms $(F)$ is defined by $(Ff)(y):= \frac1{\sqrt{2\pi}}\int_{-\infty}^\infty e^{-ixy}f(x) dx=\frac1{\sqrt{2\pi}}\int_{-\infty}^\infty (\cos xy-i\sin xy)f(x) dx.$ The Fourier cosine and Fourier sine transforms of the  function $f$, denoted by  $(F_c)$ and  $(F_s)$ respectively, are defined by integral formulas 
$(F_cf)(y):= \sqrt{\frac{2}{\pi}}\int_{\R_+} \cos (xy) f(x) dx,\ y>0,$ and $ (F_sf)(y):= \sqrt{\frac{2}{\pi}}\int_{\R_+} \sin (xy) f(x) dx,\ y>0 
$. If $f$ is an even function, the transform $(F_c)$ is the transform $(F)$ and if $f$ is an odd function, $(F_s)$ is the transform $(F)$. We define the cosine and sine of Fourier transforms in the mean-square convergence sense, namely
\begin{equation}\label{eq2.6}
(F_{\left\{\substack{c\\s}\right\}}f)(y):= \sqrt{\frac{2}{\pi}}\lim\limits_{N\to \infty}\int_{\frac1N}^N \left\{\substack{\cos xy\\ \sin xy}\right\} f(x) dx, y>0
\end{equation}
and Plancherel’s theorem in \cite{Sneddon1972} said that $F_c ,F_s : L_2 (\R_+) \leftrightarrow L_2 (\R_+)$ are
automorphism mappings with Parseval’s equalities $\|F_{\left\{{ }_s^c\right\}} f\|_{L_2(\mathbb{R}_{+} )}=\|f\|_{L_2(\mathbb{R}_{+} )}$  \cite{Titchmarsh1986}.
\begin{definition}[Sneddon's convolution \cite{Sneddon1972}] 
	The convolution for two functions $f$, $g$, related to the Fourier sine-cosine integral transforms, 
	denoted by $(\underset{1}{\ast})$, is
	 defined by
	\begin{equation}\label{eq2.9}
	\big(f\underset{1}{\ast}g\big)(x) := \frac{1}{\sqrt{2\pi}}\int_0^\infty f(x)\big[g(|x-u|)- g(x+u)\big] du,\ x>0.
	\end{equation}
	If $f,g$ are $L_1$-Lebesgue integrable functions over $\R_+$, then $\big(f\underset{1}{\ast} g\big)\in L_1(\mathbb{R}_+)$ and obtain the following factorization identity
	\begin{equation}\label{eq2.11}
	F_{s}\big(f\underset{1}{\ast} g\big)(y) = \big(F_{s}f\big)(y)(F_cg)(y),\ y>0,
	\end{equation} 
\end{definition}
\noindent Within the framework of this article, we shall make frequent use of the weighted Lebesgue spaces
$$L_p(\mathbb{R}_+,\rho):=\bigg\{f(x),\  \text{which is defined in}\ \mathbb{R_+}\ \text{such that}\ \int_0^\infty|f(x)|^p\rho(x)dx<\infty\bigg\},$$
with respect to a  positive measure $\rho(x)dx$ equipped with the norm for which $\lv f \lv_{L_p (\R_+ , \rho)}=\big\{ \int_0^{\infty} |f(x)|^p \rho(x)dx \big\}^{\frac{1}{p}},$ where $1\leq p < \infty$. If the weighted function $\rho(x) \equiv 1$ then $L_p (\R_+ ,\rho)\equiv L_p (\R_+)$.  Therefore $$L_p(\mathbb{R}_+)\subseteq L_p(\mathbb{R}_+,\rho).$$ Following \cite{Yakubovich1996index}, in this paper, the KL-transform is denoted by $(K)$ and  defined by $K[f](y):=\int_0^\infty K_{iy}(x)f(x) dx,$ $y>0,$
where $K_{iy}(x)$ is the Bessel-type function \cite{Prudnikov1986Marichev} and given by the following formula
$
K_{iy}(x)=\int_0^\infty e^{-x\cosh t}\cos( yt) dt.
$
Some results on establishing estimates for Bessel-type functions  $K_{iy}(x)$ can be found in \cite{Yaku94Luchko} which indicate that: $\forall \beta \in (0,1]$ then 
\begin{equation}\label{eq2.5}
\int_0^\infty K_{iy}(x)\cos(ty)dy=\frac{\pi}{2}e^{-x\cos t},\ \ \textup{and}\ \ |K_{iy}(x)| \leq K_0(x)\leq K_0(\beta x),\ \ \textup{and}\  |K_{iy}(x)|\leq e^{-|y|\arccos \beta}K_0(\beta y).
\end{equation}
 Here $K_0(.)$ is a Bessel-type function with the zero-exponent and defined by $K_0(x):=\int_0^\infty e^{-x\cosh t} dt, x>0$.
Based on the results in \cite{Yakubovich1996index}, it can be affirmed that, the KL-transform is an isometric isomorphism mapping  from $L_2(\mathbb{R}_{+}, x dx) \to L_2(\mathbb{R}_{+}, x \sinh \pi xdx)$, where integral $\int_{0}^{\infty} K_{iy}(x) f(x) dx$  does not exist in Lebesgue's sense, and therefore we understand it in the form 
\begin{equation}\label{eq2.7}
K[f](y):= \lim_{N\to \infty}\int_{\frac1N}^\infty K_{iy}(x) f(x) dx.
\end{equation}
 Following \cite{yakubovich2010Britvina}, we mention
$
L_p^{\alpha,\beta}(\mathbb{R}_+):=\big\{f(x),$ which is defined in $\mathbb{R}_+$ such that $\int_0^\infty|f(x)|^pK_0(\beta x)x^\alpha dx<\infty\big\}.
$ This is called a two-parameter family of Lebesgue spaces
and equipped
$
\|f\|_{L_p^{\alpha,\beta}(\mathbb{R}_+)}=\left\{\int_0^\infty|f(x)|^pK_0(\beta x)x^\alpha dx\right\}^{\frac1p},$ where $1\leq p<\infty,
$
and $\beta\in(0,1]$.  Notice that  $L_p(\mathbb{R}_+)\subset L_p^{\alpha,\beta}(\mathbb{R}_+)$, $p\geq 1$, $\alpha \in \R$ and $0<\beta\leq 1$.

\begin{definition}[Yakubovich-Britvina's convolution \cite{yakubovich2010Britvina}]\label{DN3}
	 The generalized convolution of  two functions $f$ and $g$, related to the Fourier sine-Kontorovich-Lebedev transform, denoted by $(\underset{2}{\ast})$, was  by Yakubovich-Britvina as follows
		\begin{equation}\label{eq2.12}
	\big(f \underset{2}{\ast} g\big) := \frac12\int_{\mathbb{R}^2} \big[e^{-u\cosh (v-x)} - e^{-u\cosh(v+u)}\big] f(v)g(u) du dv,\ x>0.
	\end{equation}
\noindent Furthermore, Theorem 5, p.266 in \cite{yakubovich2010Britvina} asserts that: If $f\in L_1(\mathbb{R}_+)$ and $g\in L_1^{0,\beta}(\mathbb{R}_+)$, $\forall \beta \in (0,1]$ then $\big(f \underset{2}{\ast} g\big)$  belongs to $ L_1(\mathbb{R}_+)$ and satisfies the following factorization identity
\begin{equation}\label{eq2.13}
F_s\big(f \underset{2}{\ast} g\big)(y) = (F_sf)(y)K[g](y),\ y>0.
\end{equation}	\end{definition}

\subsection{Organization}
This paper is organized into four sections. As stated in the Abstract, the primary contributions of this paper are presented in Section~\ref{sec2} and Section~\ref{sec3}. In the final part, Section~\ref{sec4} addresses the applicability of the polyconvolution operator to the solvability of a particular class of integral equations of the Toeplitz–Hankel type. This investigation is carried out by employing the results established in the preceding sections in conjunction with the Wiener–Lévy theorem in \cite{NaimarkMA1972} and associated convolution inequalities. Integral inequalities are fundamental tools for analyzing the qualitative and quantitative properties of integral transforms and differential equation solutions. In particular, convolution inequalities are essential and in fact, indispensable, as numerous integral transforms and ODE/integral equation solutions are expressed in terms of convolutions \cite{Davies2002integral}.

\section{Composition structure and mapping properties of polyconvolution related to KL-transform}\label{sec2}
\subsection{Definition, existence on $L_1$ space and  factorization identity}

Now, we apply Definition \ref{dn1} in the context where four integral transforms will be acting in the following schematic order $\overset{\gamma}{\underset{\mathcal{K}_4,\mathcal{K}_1, \mathcal{K}_2,\mathcal{K}_3}{\ast}}: \prod_{j=1}^3 U_j(X_j)\to \mathcal{V}(Y)$, such that satisfies the factorization identity
$
\mathcal{K}_4\big[\overset{\gamma}{\underset{\mathcal{K}_4,\mathcal{K}_1,\mathcal{K}_2,\mathcal{K}_3}{\ast}}(f,g,h)\big](y)=\gamma(y)(\mathcal{K}_1 f)(y)(\mathcal{K}_2 g)(y)(\mathcal{K}_3 h)(y).
$
Hence, by the factorization identity, we infer that
$
\overset{\gamma}{\underset{\mathcal{K}_4,\mathcal{K}_1,\mathcal{K}_2,\mathcal{K}_3}{\ast}}(f,g,h)(x)=\mathcal{K}_4^{-1}\big[\gamma(y)(\mathcal{K}_1f)(y)(\mathcal{K}_2 g)(y)(\mathcal{K}_3 h)(y)\big](x).
$ Here always assume that $(\mathcal{K}_4)$ has an inverse transform $(\mathcal{K}_4^{-1})$. We propose polyconvolution by choosing $\mathcal{V}(Y)$ as the algebra of all Lebesgue measurable functions over $\mathbb{R_+}$; and $U_j(X_j)$ as $L_1 (\R_+)$; and the weight function $\gamma(y)$ is homogenous equal $1$ for every $y$ in half-axis $\mathbb{R}_+$. The integral transforms are selected in order as follows: $(\mathcal{K}_4)\equiv (\mathcal{K}_1)$ is the Fourier sine transform; $(\mathcal{K}_2)$ is the  Fourier cosine transform; and $(\mathcal{K}_3)$ is the Kontorovich-Lebedev transform.
\begin{definition}[\textsc{Main object}]\label{DN4}
The polyconvolution operator of three functions $f$, $g$ and $h$ generated by integral transforms including Fourier sine, Fourier cosine, and Kontorovich-Lebedev is denoted by $\underset{F_s,F_c,K}{\ast}(f,g,h)$ and defined by the formula	
\begin{equation}\label{eq3.1}
\underset{F_s,F_c,K}{\ast}(f,g,h)(x):=\frac1{2\sqrt{2\pi}}\int_{\mathbb{R}_+^3}\Phi(x,u,v,w)f(u)g(v)h(w)du dv dw,\ x>0,
\end{equation}
where kernel function $\Phi$ is represented by
\begin{equation}\label{eq3.2}
\Phi(x,u,v,w)=e^{-w\cosh(x-u+v)} + e^{-w\cosh(x-u-v)}-e^{-w\cosh(x+u+v)}-e^{-w\cosh(x+u-v)}.
\end{equation}
\end{definition}
\noindent Definition \ref{DN4} is a successor and extension of Definition \ref{DN3} in terms of operator structure. First, we have the estimate of the kernel function.
\begin{lemma}
Setting $\mathcal{I}_1=\int_0^\infty |\Phi(x,u,v,w)|dx$; $\mathcal{I}_2=\int_0^\infty |\Phi(x,u,v,w)|du$; $\mathcal{I}_3=\int_0^\infty|\Phi(x,u,v,w)|dv$; and\\ $\mathcal{I}_4=\int_0^\infty|\Phi(x,u,v,w)|dw.$ Then, we obtain 
\begin{equation}
\label{eq3.3} \mathcal{I}_i\leq 4K_0(w),\ \forall  i=\overline{1,3}\ \  \textup{and}\ \  
\mathcal{I}_4\leq \frac{4}{\cosh t}=4\,\mathrm{sech}\,t,
\end{equation}
where $\cosh t=\min\big\{\cosh(x-u+v),\cosh(x-u-v),\cosh(x+u+v),\cosh(x+u-v)\big\}$.
\end{lemma}
\begin{proof}
	From \eqref{eq3.2} and performing the change of variables, we have
\begin{align*}
&\mathcal{I}_1\leq \int_0^\infty e^{-w\cosh(x-u+v)}dx + \int_0^\infty e^{-w\cosh(x-u-v)}dx + \int_0^\infty e^{-w\cosh(x+u+v)}dx + \int_0^\infty e^{-w\cosh(x+u-v)}dx\\
&= \int_{-(u-v)}^\infty e^{-w\cosh t}dt + \int_{-(u+v)}^\infty e^{-w\cosh t} dt + \int_{u+v}^\infty e^{-w\cosh t}dt + \int_{u-v}^\infty e^{-w\cosh t}dt
\leq 2\int_{-\infty}^\infty e^{-w\cosh t}dt=4K_0(w).
\end{align*}
Similar to the above method,  we obtain $\mathcal{I}_2,\mathcal{I}_3<4K_0(w), \forall w>0$.  Meanwhile
\begin{align*}
\mathcal{I}_4&\leq\int_0^\infty e^{-w\cosh(x-u+v)} dw + \int_0^\infty e^{-w\cosh(x-u-v)}dw
 + \int_0^\infty e^{-w\cosh(x+u+v)}dw + \int_0^\infty e^{-w\cosh(x+u-v)}dw\\
&=\frac{1}{\cosh(x-u+v)} + \frac{1}{\cosh(x-u-v)} + \frac{1}{\cosh(x+u+v)} + \frac{1}{\cosh(x+u-v)}
\leq \frac{4}{\cosh t}=4\mathrm{sech}\,t,
\end{align*}
\end{proof}
\begin{theorem}\label{thm3.1}
Suppose that $f,g$ are arbitrary functions in $L_1 (\R_+ )$, and $h\in L_1^{0,\beta}(\mathbb{R}_+)$, we have two assertion as follows:\\
\textsc{A)} For  $\beta \in (0,1]$, then the polyconvolution \eqref{eq3.1} is well-defined  for all $x >0$ as a continuous function and belongs to $L_1(\mathbb{R}_+)$ and we obtain the following estimate 
\begin{equation}\label{eq3.7}
\big\|\underset{F_s,F_c,K}{\ast} (f,g,h)\big\|_{L_1(\mathbb{R}_+)} \leq \sqrt{\frac{2}{\pi}}\|f\|_{L_1(\mathbb{R}_+)} \|g\|_{L_1(\mathbb{R}_+)} \|h\|_{L_1^{0,\beta}(\mathbb{R}_+)}.
\end{equation}
\textsc{B)} In case $0<\beta<1$, for all $x>0$, then polyconvolution \eqref{eq3.1} satisfies the generalized Parseval-type identity 
\begin{equation}\label{eq3.5}
\underset{F_s,F_c,K}{\ast}(f,g,h)(x)=\sqrt{\frac{2}{\pi}}\int_0^\infty (F_sf)(y)(F_cg)(y)K[h](y)\sin(xy) dy,\ x>0,
\end{equation}
and the following factorization property is valid
\begin{equation}\label{eq3.6}
F_s\big(\underset{F_s,F_c,K}{\ast} (f,g,h)\big)(y) = (F_sf)(y)(F_cg)(y)K[h](y), y>0.
\end{equation}	
Moreover $\underset{F_s,F_c,K}{\ast} (f,g,h)\in C_0(\mathbb{R}_+)$, where $C_0(\mathbb{R}_+)$ is the space of bounded continuous functions vanishing at infinity.
\end{theorem}
\begin{proof} To prove assertion \textsc{A)}, we need to show $
	\int_0^\infty \big|\underset{F_s,F_c,K}{\ast}(f,g,h)(x)\big|dx$ is finite.
	Indeed, from the formulas \eqref{eq3.1}, \eqref{eq3.2},  we obtain
	$
	\int_0^\infty\big|\underset{F_s,F_c,K}{\ast}(f,g,h)(x)\big|dx\leq \frac{1}{2\sqrt{2\pi}}\int_{\mathbb{R}_+^4}|\Phi(x,u,v,w)||f(u)| |g(v)| |h(w)| du dv dw dx
	$. By estimate of kernel function $\Phi$ in \eqref{eq3.3}, coupling with Fubini’s theorem, we obtain 
	$$
	\begin{aligned}
	\int_0^\infty\big|\underset{F_s,F_c,K}{\ast}(f,g,h)(x)\big|dx &\leq \frac{1}{2\sqrt{2\pi}}\left(\int_0^\infty |f(u)|du\right)\left(\int_0^\infty|g(v)|dv\right)\int_0^\infty |h(w)|\left(\int_0^\infty|\Phi(x,u,v,w)|dx\right)dw\\
	&\leq \frac{2}{\sqrt{2\pi}}\left(\int_0^\infty |f(u)|du\right)\left(\int_0^\infty |g(v)|dv\right) \left(\int_0^\infty K_0(w)|h(w)|dw\right)\\
	&\leq \sqrt{\frac{2}{\pi}}\left(\int_0^\infty |f(u)|du\right)\left(\int_0^\infty|g(v)|dv\right)\left(\int_0^\infty K_0(\beta w)|h(w)|dw\right)\\
	&=\sqrt{\frac{2}{\pi}}\|f\|_{L_1(\mathbb{R}_+)}\|g\|_{L_1(\mathbb{R}_+)}\|h\|_{L_1^{0,\beta}(\mathbb{R}_+)}.
	\end{aligned}$$
Since 	$\|h\|_{L_1^{0,\beta}(\mathbb{R}_+)}= \int_0^\infty|h(x)| K_0(\beta x)dx$ is finite, then  $
\int_0^\infty  \big|\underset{F_s,F_c,K}{\ast}(f,g,h)(x)\big|dx$  is finite
for almost all $x>0$, and it implies that polyconvolution \eqref{eq3.1}  belongs to $L_1 (\R_+)$ and we derive the estimate \eqref{eq3.7}.\\
To prove \textsc{B)}, first, 
coupling \eqref{eq3.2} and \eqref{eq2.5} we have
	\begin{align*}
	\Phi(x,u,v,w)&=\frac{2}{\pi}\int_0^\infty K_{iy}(w)\big[\cos y(x-u+v)+\cos y(x-u-v)
-\cos y(x+u+v) - \cos y(x+u-v)\big]dy\\
	&=\frac{8}{\pi}\int_0^\infty \sin(yu)\cos(yv)\sin(xy)K_{iy}(w) dy.
	\end{align*}
	Moreover, from \eqref{eq3.1} we get
	\begin{equation}\label{eq3.9}
	\underset{F_s,F_c,K}{\ast}(f,g,h)(x)=\frac{1}{2\sqrt{2\pi}}\frac{8}{\pi}\int_{\mathbb{R}_+^4} \sin(yu)\cos(yv)\sin(xy)K_{iy}(w) f(u)g(v)k(w) du dv dw dy, x>0.
	\end{equation}
For any $f,g$ be $L_1$-Lebesgue integrable functions over $\R_+$ and $h\in L_1^{0,\beta}(\mathbb{R}_+)$, using \eqref{eq2.5} $\forall \beta \in (0,1)$, we infer that
	\begin{align*}
	\big|\underset{F_s,F_c,K}{\ast}(f,g,h)(x)\big|&\leq \frac{4}{\pi\sqrt{2\pi}}\int_{\mathbb{R}_+^4} e^{-y \arccos\beta}K_0(\beta w)|f(u)| |g(v)| |h(w)| du dv dw dy\\
	&=\frac{4}{\pi\sqrt{2\pi}}\left(\int_0^\infty|f(u)|du\right)\left(\int_0^\infty |g(v)|dv\right)\left(\int_0^\infty K_0(\beta w)|h(w)| dw\right)\left(\int_0^\infty e^{-y\arccos \beta}dy\right)\\
	&=\frac{4}{\pi\sqrt{2\pi}}\|f\|_{L_1(\mathbb{R}_+)} \|g\|_{L_1(\mathbb{R}_+)}\|h\|_{L_1^{0,\beta}(\mathbb{R}_+)} \frac{1}{\arccos \beta}<\infty.
	\end{align*}
	This implies that the integration \eqref{eq3.9} is absolutely convergent, by using Funini's theorem, we have
$$	\begin{aligned}
	&\underset{F_s,F_c,K}{\ast}(f,g,h)(x)=\frac{4}{\pi\sqrt{2\pi}}\int_0^\infty\left\{\int_0^\infty f(u)\sin (yu)\left(\int_0^\infty g(v)\cos (yv)
	\left(\int_0^\infty K_{iy}(w)h(w)\right)dv\right)du\right\}\sin (xy) dy,\\
	&=\sqrt{\frac{2}{\pi}}\int_0^\infty\left\{\left(\sqrt{\frac{2}{\pi}}\int_0^\infty f(u)\sin (yu) du\right)\left(\sqrt{\frac{2}{\pi}}\int_0^\infty g(v) \cos (yv) dv\right)
	\left(\int_0^\infty K_{iy}(w)h(w)dw\right)\right\}\sin (xy) dy\\
	&=\sqrt{\frac{2}{\pi}}\int_0^\infty (F_sf)(y)(F_cg)(y)(K[h])(y)\sin (xy) dy.
	\end{aligned}$$
	Therefore, we deduce that the Parseval identity \eqref{eq3.5} holds.  Applying the Fourier sine transform $(F_s)$ on both sides of the Parseval identity, we derive the factorization equality \eqref{eq3.6}. 
	  Riemann-Lebesgue's theorem in \cite{Sogge1993fourier} states that \textquotedblleft If $f\in L_1(\R^n)$, then $(Ff)(y) \to 0$ as  $|y|\to\infty$, and, hence $(Ff)(y)\in C_0 (\R^n)$\textquotedblright. This  assertion is still true for  Fourier sine transform  $(F_s)$ and  Fourier cosine transform $(F_c)$  on $\R_+$  \cite{Titchmarsh1986}, implying that if $f \in L_1(\mathbb{R}_+)$, then $(F_s f)\in C_0(\mathbb{R}_+)$ and $|(F_sf)(y)|\leq \|f\|_{L_1(\mathbb{R}_+)}$, similarly if $g \in L_1(\mathbb{R}_+)$, then $(F_c g)\in C_0(\mathbb{R}_+)$ and $|(F_c g)(y)|\leq \|f\|_{L_1(\mathbb{R}_+)}$.
Following \eqref{eq2.5} we deduce 
	$|K[h](y)|\leq \int_0^\infty K_{iy}(x)|h(x)| dx \leq \int_0^\infty e^{-y\arccos\beta}K_0(\beta x)|h(x)|dx
	=e^{-y\arccos\beta}\|h\|_{L_1^{0,\beta}(\mathbb{R}_+)}$ is finite.
	Therefore, $(F_sf)(y)(F_cg)(y)K[h](y)$ is  a bounded function on  ${\R}_+$. Letting $x$ tend to $\infty$ in Parseval's identity \eqref{eq3.5}, we obtain the final conclusion of this theorem obviously by virtue of Riemann-Lebesgue's theorem.
\end{proof}
\begin{remark}[Connection with previously known convolutions]\label{rmk1}\textup{Observing the factorization identities \eqref{eq2.11}, \eqref{eq2.13} and \eqref{eq3.6}, we obtain the correlation between polyconvolution  \eqref{eq3.1} with Sneddon's convolution \eqref{eq2.9} and Yakubovich-Britvina's convolution \eqref{eq2.12} as follows:
		Let $f,g$ be $L_1$-Lebesgue integrable functions over $\R_+$ and $h\in L_1^{0,\beta}(\mathbb{R}_+)$, $\forall \beta \in (0,1]$, then
			$
			\underset{F_s,F_c,K}{\ast}(f,g,h)(x)= \big[\big(f\underset{1}{\ast}g\big)\underset{2}{\ast}h\big](x)$ belonging to $L_1(\mathbb{R}_+).
			$ Therefore, the composition structure of operator \eqref{eq3.1} offers more favorable properties than the structure of generalized convolutions introduced in \cite{TuanMMA2024,trigonometric2024}.}
\end{remark}
\begin{theorem}\label{thm3.2}
	Let $f,g$ be the quadratically integrable functions over $\mathbb{R}_+$ and $h\in L_2^{0,\beta}(\mathbb{R}_+)$, $0<\beta<1$, then polyconvolution operator \eqref{eq3.1} satisfies Parseval's identity \eqref{eq3.5} and factorization identity \eqref{eq3.6} (the integrals here are understood in the sense of Cauchy principal value), where $(F_s)$, $(F_c)$ transforms in $L_2$ are defined by \eqref{eq2.6} and KL-transform in $L_2$ is defined by \eqref{eq2.7}.
\end{theorem}
\begin{proof}
	Arguing similarly as in the proof of Theorem \ref{thm3.1}, we obtain Parseval's identity \eqref{eq3.5}. Since the transforms $F_s, F_c: L_2 (\R_+)\leftrightarrow L_2 (\R_+)$ are isometric isomorphism mappings \cite{Titchmarsh1986}, it follows that
	if $f,g\in L_2(\mathbb{R}_+)$, then  $F_sf, F_cg$ belong to $L_2(\mathbb{R}_+)$. By estimate \eqref{eq2.5}, with $0<\beta<1$, we infer that
	\begin{align*}
	|K[h](y)|&\leq \int_0^\infty |K_{iy}(x)| |h(x)| dx\leq \int_0^\infty e^{-|y|\arccos\beta}K_0(\beta x)|h(x)| dx,\\
	&\leq \left\{\int_0^\infty e^{-|y|\arccos\beta}K_0(\beta x)dx\right\}^{\frac12}\left\{\int_0^\infty e^{-|y|\arccos\beta}K_0(\beta x)|h(x)|^2dx\right\}^{\frac12}\\
	&=e^{-|y|\arccos \beta}\left\{\int_0^\infty K_0(\beta x)dx\right\}^{\frac12}\left\{\int_0^\infty K_0(\beta x)|h(x)|^2dx\right\}^{\frac12}\\
	&\leq \textup{Const.}\|h\|_{L_2^{0,\beta}(\mathbb{R}_+)}<\infty,\ \ \textup{where}\ \ \textup{Const}:=e^{-|y|\arccos\beta}\left(\int_0^\beta K_0(\beta x)dx\right)^{\frac12}.
	\end{align*}
	Therefore, $(F_sf)(y)(F_cg)(y)K[h](y)$ belongs to $L_2(\mathbb{R}_+)$. By applying the Fourier-sine transform to both sides of the Parseval identity \eqref{eq3.5} we get \eqref{eq3.6}.
\end{proof}

\subsection{One-dimensional Watson type transform}
According to \cite{WatsonGN1933}, Watson  showed that the Mellin convolution type transforms $g(x)=[K f](x)=\int_{0}^{\infty} k(x y) f(y) d y,$
such that their inverses have a similar form
$f(x)=[\hat{K} g](x)=\int_{0}^{\infty} \hat{k}(x y) g(y) d y,$ where the kernel $ \hat{k}(x)$ is called the conjugate kernel.  Let us consider the above transform  from the point of view of general Fourier transforms \cite{Titchmarsh1986}, it is well-known that the  transform  $g(x)=\frac{\partial}{\partial x}\int_{0}^{\infty}k_1 (xy)f(y)\frac{dy}{y}$ is an automorphism over $L_2 (\R_+ ,dx)\equiv L_2 (0, \infty)$ and has the symmetric inversion formula if and only if the function $k_1 (x)$ satisfies conditions
$$\int_{0}^{\infty} k_1 (ax) k_1 (bx) x^{-2}dx = \min (a,b)\ \ \text{and} \  \
\int_{0}^{\infty} k_1 (ax) \overline{k_1 (bx)}x^{-2}dx =\min (a,b),$$
with all positive numbers $a$ and $b$. The function $k_1 (x)$ satisfying the above-mentioned conditions is called as a one-dimensional Watson kernel \cite{Titchmarsh1986}.  The first equality in Watson's conditions is the classical definition of the Watson kernel \cite{WatsonGN1933}, which leads to a unitary Watson transform when $k_1 (x)$ is a real function. In this sense, we can study the Watson-type integral transform for other convolutions as follows $f\mapsto g = D(f*k)$, where $D$ is an arbitrary differential operator and $k$ is the known kernel. In this part, we study the Watson-type transform for polyconvolution \eqref{eq3.1} by fixing a function and letting the remaining functions vary in certain function spaces. Indeed, we fix $g=g_0\in L_2(\mathbb{R}_+)$, and $h=h_0\in L_2^{0,\beta}(\mathbb{R}_+)$, $\beta\in(0,1)$ and  consider the operator $(T_{g_0,h_0}): L_2(\mathbb{R}_+)\to~L_2(\mathbb{R}_+)$ such that $f\mapsto\varphi:=D\big(\underset{F_s,F_c,K}{\ast}(f,g_0,h_0)\big)$ where, $D:=\big(I- d^2 / dx^2\big)$, with $I$ is a unitary operator and 
\begin{equation}\label{eq6.1}
\begin{aligned}
\varphi(x)&=\left(T_{g_0,h_0}f\right)(x)=\left(1-d^2 / dx^2\right)\big(\underset{F_s,F_c,K}{\ast}(f,g_0,h_0)\big)(x)\\
&=\frac{1}{2\sqrt{2\pi}}\left(1-d^2 / dx^2\right)\int_{\mathbb{R}_+^3} \Phi(x,u,v,w)f(u)g_0(v)h_0(w) du dv dw,\ x>0.
\end{aligned}
\end{equation}
The main theorem of this part is stated as follows.
\begin{theorem}\label{thm6.1}
	Let $g_0 \in L_2(\mathbb{R}_+)$, and  $h_0$ belongs to  $L_2^{0,\beta}(\mathbb{R}_+)$, $\beta\in(0,1)$  satisfy 
\begin{equation}\label{eq6.2}
\big|(F_cg_0)(y)K[h_0](y)\big| = \frac{1}{1+y^2},\ y>0.
\end{equation}	
	Then the condition \eqref{eq6.2} is the necessary and sufficient condition for operator $(T_{g_0,h_0})$ to become an isometric isomorphism mapping on $L_2 (\R_+)$. Moreover, the inverse operator of $(T_{g_0,h_0})$ has a symmetric form and is represented by 
	\begin{equation}\label{eq6.3}
	f(x)=\left(T_{g_0,h_0}^{-1}\varphi\right)(x)\equiv \left(T_{\bar{g}_0,\bar{h}_0}\varphi\right)(x)=\left(1-d^2 / dx^2\right)\big(\underset{F_s,F_c,K}{\ast}(\varphi,\bar{g}_0,\bar{h}_0)\big)(x),
	\end{equation}
	where, $\bar{g}_0$, $\bar{h}_0$ are conjugate functions of $g_0$, $h_0$ respectively.
\end{theorem}
\begin{proof}
\textsc{Necessary condition:} Theorem 68 in \cite{Titchmarsh1986} states that \textquotedblleft If $f(y)$, $yf(y)$ and $y^2f(y)\in L_2(\mathbb{R}_+)$ then $(Ff)(x)$, $\frac{d}{dx}(Ff)(x)$, and  $\frac{d^2}{dx^2}(Ff)(x)\in L_2(\mathbb{R}_+)$\textquotedblright. This  assertion is still true for  Fourier sine transform  $(F_s)$. Meaning that, if $f(y)$, $yf(y)$ and $y^2f(y)\in L_2(\mathbb{R}_+)$ then $(F_sf)(x)$, $\frac{d}{dx}(F_sf)(x)$ and $\frac{d^2}{dx^2}(F_sf)(x)\in L_2(\mathbb{R}_+)$. Applying \eqref{eq2.6}, we have 	\begin{align*}
\frac{d}{dx}(F_sf)(x)&=\sqrt{\frac{2}{\pi}}~\lim_{N\to \infty} \int_{\frac{1}{N}}^N f(y)\frac{d}{dx}\sin(xy) dy
\\&= \sqrt{\frac{2}{\pi}}~\lim_{N\to \infty} \int_{\frac1N}^N yf(y)\cos(xy)dy=F_c(yf(y))(x)\in L_2(\mathbb{R}_+).\\
\text{Likewise,}\quad \frac{d^2}{dx^2}(F_sf)(y)&=\frac{d}{dx}\left(\frac{d}{dx}(F_sf)(y)\right)(x)
=\frac{d}{dx}\big[F_c(yf(y))\big](x)=\sqrt{\frac{2}{\pi}}\lim_{N\to\infty}\int_{\frac1N}^N yf(y) \frac{d}{dx}\cos(xy)dy\\
&=\sqrt{\frac{2}{\pi}}\lim_{N\to \infty}\int_{\frac1N}^N -y^2f(y)\sin(xy)dy=F_s(-y^2f(y))(x)\in L_2(\mathbb{R}_+).
\end{align*}
Therefore, we claim that
\begin{equation}\label{eq6.4}
\left(1-d^2 / dx^2\right)(F_sf)(y)=F_s((1+y^2)f(y))(x)\in L_2(\mathbb{R}_+).
\end{equation}	On the other hand, following condition \eqref{eq6.2}, we deduce $(1+y^2)(F_cg_0)(y)K[h_0](y)$ is a bounded function, implying that $(1+y^2)(F_cg_0)(y)K[h_0](y) (F_sf)(y)$ belongs to $L_2(\mathbb{R}_+)$ holds $\forall f\in L_2(\mathbb{R}_+).$ Using the formulas \eqref{eq6.1}, \eqref{eq6.4} and Parseval's identity \eqref{eq3.5} we infer that
\begin{equation}\label{eq6.5}
\begin{aligned}
\varphi(x)=\left(T_{g_0,h_0}f\right)(x)&=\left(1-d^2 / dx^2\right)\big(\underset{F_s,F_c,K}{\ast}(f,g_0,h_0)\big)(x)\\
&=F_s\left\{(1+y^2)(F_sf)(y)(F_cg_0)(y)K[h_0](y)\right\}\in L_2(\mathbb{R}_+).
\end{aligned}
\end{equation}
Based on Theorem 9.13 in \cite{WRudin1987}, for any $f$ be a function belonging to the $L_2 (\R_+)$, then $\lv Ff\lv_{L_2 (\R_+)}=\lv f\lv_{L_2 (\R_+)}$.  This is still true for the $(F_s)$ transform which implies that $\lv F_s f\lv_{L_2 (\R_+)}=\|f\|_{L_2 (\R_+)}, \forall f\in L_2 (\R_+)$, under the condition \eqref{eq6.2}, we obtain
\begin{align*}
\|\varphi\|_{L_2(\mathbb{R}_+)}=\left\|T_{g_0,h_0} f\right\|_{L_2(\mathbb{R}_+)} &=\left\|F_s\left\{(1+y^2)(F_sf)(y)(F_cg_0)(y)K[h_0](y)\right\}\right\|_{L_2(\mathbb{R}_+)}\\
&=(1+y^2)\left|(F_cg_0)(y)K[h_0](y)\right| \|F_sf\|_{L_2(\mathbb{R}_+)}=\|F_sf\|_{L_2(\mathbb{R}_+)}=\|f\|_{L_2(\mathbb{R}_+)}.
\end{align*}
This means $(T_{g_0,h_0})$ is an isometric isomorphism  transformation (unitary) in $L_2(\R_+)$. Now, we need to show that the inverse operator of $T_{g_0,h_0}$ has a symmetric-type form. Indeed, in $L_2$, by applying the unitary property of mapping $F_s: L_2 (\R_+)\leftrightarrow L_2 (\R_+)$ for equality \eqref{eq6.5}, we deduce that
\begin{equation}\label{eq6.6}
(F_s\varphi)(y)=F_s\left(T_{g_0,h_0}f\right)(y)=(1+y^2)(F_sf)(y)(F_cg_0)(y)K[h_0](y)\in L_2(\mathbb{R}_+).
\end{equation}
On the other hand, since $\{\overline{(F_cg_0)}(y)\cdot \overline{K[h_0]}(y)=(F_c\bar{g}_0)(y)K[\bar{h}_0](y)\}$.
Multiplying both sides of \eqref{eq6.6} by expression $(F_c\bar{g}_0)(y)K[\bar{h}_0](y)\}$ and applying \eqref{eq6.2}, we get
$
(F_s\varphi)(y)|(F_c\bar{g}_0)(y) K[\bar{h}_0](y)|=\frac{1}{1+y^2}(F_sf)(y)
$
and
\begin{equation}\label{eq6.7}
(F_sf)(y)=(1+y^2)|(F_c\bar{g}_0)(y)K[\bar{h}_0](y)|(F_s\varphi)(y).
\end{equation}
By making the same argument as above, we get $F_s \varphi \in L_2 (\R_+)$ and  $(1+y^2)(F_c\bar{g}_0)(y)K[\bar{h}_0](y)$  is a bounded function, implying
that $(1+y^2)(F_c\bar{g}_0)(y)K[\bar{h}_0](y)\allowbreak(F_s\varphi)(y) $ belongs to $L_2(\mathbb{R}_+)$, $\forall \varphi\in L_2(\mathbb{R}_+)$. Combining \eqref{eq6.4}, \eqref{eq6.7} and Parseval's identity \eqref{eq3.5}, we obtain
\begin{align*}
f(x)=F_s\left\{(1+y^2)(F_s\varphi)(y) (F_c\bar{g}_0)(y) K[\bar{h}_0](y)\right\}(x)
&=\left(1-d^2 / dx^2\right)F_s\left\{(F_s\varphi)(y) (F_c\bar{g}_0)(y) K[\bar{h}_0](y)\right\}(x)\\
&=\left(1-d^2 / dx^2\right)\big(\underset{F_s,F_c,K}{\ast}(\varphi,\bar{g}_0,\bar{h}_0)\big)(x)=\left(T_{\bar{g}_0,\bar{h}_0}\varphi\right)(x).
\end{align*}
\textsc{Sufficient condition:} Suppose that $(T_{g_0,h_0})$ is a unitary operator on $L_2(\mathbb{R}_+)$ and has the inverse operator $(T_{\bar{g}_0,\bar{h}_0})$. We need to prove that the functions $g_0,h_0$ must satisfy the condition \eqref{eq6.2}. Indeed, since $T_{g_0,h_0}$ has the unitary property on $L_2(\R_+)$, then for any functions $f$ belonging to $L_2 (\R_+)$, we obtain 
\begin{align}\label{eq6.8}
\|\varphi\|_{L_2(\mathbb{R}_+)}=\left\| T_{g_0,h_0}f\right\|_{L_2(\mathbb{R}_+)}=&\|f\|_{L_2(\mathbb{R}_+)}=\|F_sf\|_{L_2(\mathbb{R}_+)}\nonumber\\
&=\|(F_sf)(y)\{(1+y^2)(F_cg_0)(y)K[h_0](y)\}\|_{L_2(\mathbb{R}_+)}.
\end{align} This shows that there exists a multiplication operator of the form $\mathcal{N}_{\Theta}[.]$  defined by $
\mathcal{N}_{\Theta}[f](y):=\Theta(y). f(y),
$ with the function
 $\Theta(y)=(1+y^2)(F_cg_0)(y)K[h_0](y), \forall y>0$. Therefore, for any $f\in L_2(\mathbb{R}_+)$ the formula \eqref{eq6.8} can always be rewritten as
$
\|F_sf\|_{L_2(\mathbb{R}_+)} = \|\mathcal{N}_{\Theta}(F_sf)\|_{L_2(\mathbb{R}_+)}$.
This means that the  multiplication operator $\mathcal{N}_{\Theta}(.)$ is an isometric isomorphism  on the $L_2 (\R_+)$, and this happens if and only if
$
|(F_cg_0)(y)K[h_0](y)|=\frac{1}{1+y^2}.
$ Hence, both functions $g_0$ and $h_0$ must satisfy \eqref{eq6.2}. In summary, condition \eqref{eq6.2} is necessary and sufficient for $(T_{g_0,h_0})$ to ensure a unitary transformation on $L_2 (\R_+)$.

\end{proof}
\subsection{Plancherel theorem for Watson-type transform}
Plancherel's theorem was proved in \cite{plancherel1910contribution}, where they  stated that: \textquotedblleft the integral of a function's squared modulus is equal to the integral of the squared modulus of its frequency spectrum\textquotedblleft. That is, if $f(x)$ is a function on the real line, and $\widehat{f}(\xi)$ is its frequency spectrum, then
$
\int_{-\infty}^{\infty}|f(x)|^2 d x=\int_{-\infty}^{\infty}|\widehat{f}(\xi)|^2 d \xi.
$
A more precise formulation is that if a function is in both $L_p$ spaces $L_1(\mathbb{R})$ and $L_2(\mathbb{R})$, then its Fourier transform is in $L_2(\mathbb{R})$ and the Fourier transform is an isometry concerning the $L_2$ norm. This implies that the Fourier transform restricted to $L_1(\mathbb{R}) \cap L_2(\mathbb{R})$ has a unique extension to a linear isometric mapping $L_2(\mathbb{R}) \mapsto L_2(\mathbb{R})$, and is called the Plancherel transform (see \cite{WRudin1987}, Chapter 9, Theorem 9.13). This isometry is a unitary mapping. Afterward, in 1949, Lebedev proved the Plancherel theorem for the Kontorovich–Lebedev transform \cite{Lebedev1949parseval}. Being influenced by those mentioned above, in this part, we are motivated to obtain Plancherel-type theorem for two operators $(T_{g_0,h_0})$ and $(T_{\bar{g}_0,\bar{h}_0})$. We show that these operators can be convergent in $L_2$-norm generated by the sequence of integral operators $\{\varphi_N\}$, $\{f_N\}$. The last part deals with the boundedness of operator $(T_{g_0,h_0})$ from $L_p(\mathbb{R}_+)$ to $ L_q(\mathbb{R}_+)$, where 
$p$, $q$ are paired conjugate exponents with $p\in [1,2]$, via Riesz-Thorin's interpolation theorem in \cite{Duoandikoetxea01}. 
\begin{theorem}\label{thm6.2}
Suppose that $g_0\in L_1(\mathbb{R}_+) \cap L_2(\mathbb{R}_+)$, $h_0\in L_1^{0,\beta}(\mathbb{R}_+)\cap L_2^{0,\beta}(\mathbb{R}_+)$, $\beta \in (0,1)$ satisfies condition \eqref{eq6.2} such that $\Psi(x)=\left(1-d^2 / dx^2\right)\Phi(x,u,v,w)$ is a locally bounded function $\forall x\in \mathbb{R}_+$. Let f be a square-integrable function over $\R_+$ for any positive integer $N$, we set
\begin{equation}\label{eq6.9}
\varphi_N(x):=\frac{1}{2\sqrt{2\pi}}\int_{\mathbb{R}_+^3}\Psi(x)f^N(u) g_0(v) h_0(w) du dv dw,\ x>0,
\end{equation}
where $f^N:=f. \mathcal{X} _{[0,N]}$.  Then the following two statements hold:\\
\textsc{A)} $\varphi_N (x)$ belongs to the $L_2(\mathbb{R}_+)$, and the sequence of functions $\{\varphi_N (x)\}$ converges in the sense of norm to a function $\varphi(x)$ in $L_2 (\R_+)$ when $N\to\infty$ satisfies $\|\varphi\|_{L_2 (\R_+)}=\|f\|_{L_2 (\R_+)}.$ \\
\textsc{B)}  Let	$\varphi^N:=\varphi.\mathcal{X} _{[0,N]},$ and set
 $f_N(x)=\left(1-d^2 / dx^2\right)\big(\underset{F_s,F_c,K}{\ast} (\varphi^N, \bar{g}_0, \bar{h}_0)\big)(x),
$ then $\{f_N(x)\} \to f(x)$ in $L_2 (\R_+)$ when $N$ tends to $\infty$.
\end{theorem}
\begin{proof}
 Since $f^N = f. \mathcal{X}_{[0, N]},$ it implies that $f^{N} \in L_2 (\R_+)$, which $\mathcal{X}_{[0, N]}$ is the characteristic function of $f$ over finite interval $[0, N]$. By the assumptions of Theorem \ref{thm6.2}, then integration in formula \eqref{eq6.9}  is the
 absolute convergence. We change the order of integration and differentiation here, thus we obtain 
 \begin{equation}\label{eq6.11}
 \begin{aligned}
 \varphi_N(x)&=\frac{1}{2\sqrt{2\pi}}\int_{\mathbb{R}_+^3}\left(1-d^2 / dx^2\right)\Phi(x,u,v,w) f^N(u) g_0(v) h_0(w) du dv dw,\\
 &=\frac{1}{2\sqrt{2\pi}}\left(1-d^2 / dx^2\right)\int_{\mathbb{R}_+^3} \Phi(x,u,v,w) f^N(u) g_0(v) h_0(w) du dv dw\\
 &=\left(1-d^2 / dx^2\right)\big(\underset{F_s,F_c,K}{\ast}(f^N, g_0,h_0)\big)(x)=\big(T_{g_0,h_0}f^N\big)(x).
 \end{aligned}
 \end{equation}
 Due to Theorem \ref{thm6.1}, we infer that $\varphi_N\in L_2(\mathbb{R}_+)$.  Based on the setting of $f^N$, obviously that $\|f^N-f\|_{L_2(\mathbb{R}_+)}\to 0$, when $N\to\infty$, coupling with \eqref{eq6.1}	 then
 $$
(\varphi_N-\varphi)(x)=\left(1-d^2 / dx^2\right)\big(\underset{F_s,F_c,K}{\ast}(f^N-f,g_0,h_0)\big)(x)= \big(T_{g_0,h_0}(f^N-f)\big)(x).
 $$
On the other hand, since  $f^N-f\in L_2(\mathbb{R}_+)$ and the functions $g_0,h_0$ satisfy condition \eqref{eq6.2}, which implies that $\varphi_N-\varphi\in L_2(\mathbb{R}_+)$ and we obtain 
 $
 \|\varphi_N-\varphi\|_{L_2(\mathbb{R}_+)}=\left\|T_{g_0,h_0}(f^N-f)\right\|_{L_2(\mathbb{R}_+)}	=\|f^N-f\|_{L_2(\mathbb{R}_+)}.
 $
This means that $\varphi_N\to \varphi\in L_2(\mathbb{R}_+)$, when $N$ tends to $\infty$. Moreover, seeing $\varphi$ as the image of $f$ under the operator $(T_{g_0,h_0})$ in $L_2 (\R_+)$, by Theorem \ref{thm6.1}, we obtain $\lVert f\lVert_{L_2 (\R)}=\left\|T_{g_0,h_0} f\right\|_{L_2(\mathbb{R}_+)}=\| \varphi\|_{L_2 (\R_+)}$, which is unitary. So the proof of the first assertion is completed.
 
 By reasoning similarly to the above for  $f(x)=(T_{\bar{g}_0,\bar{h}_0}\varphi)(x)$, as defined by \eqref{eq6.3}, where, $\bar{g}_0$, $\bar{h}_0$ are conjugate functions of $g_0$, $h_0$ respectively, and satisfy the condition \eqref{eq6.2}, we have
 $$\begin{aligned}
 (\varphi_N-\varphi)(x)=\left(1-d^2 / dx^2\right)\big(\underset{F_s,F_c,K}{\ast}(\varphi_N-\varphi,\bar{g}_0,\bar{h}_0)\big)(x)=\big(T_{\bar{g}_0,\bar{h}_0}(\varphi^N-\varphi)\big)(x). 
 \end{aligned}
 $$
 Thus, $\varphi\in L_2\mathbb{R}_+$ then $\varphi^N:=\varphi. \mathcal{X}_{[0,N]}\in L_2(\mathbb{R}_+)$, this leads to $\varphi^N-\varphi \in L_2(\mathbb{R}_+)$. By applying Theorem \ref{thm6.1}, we have $
 f_N-f\in L_2(\mathbb{R}_+)
 $
 and
 $
 \|f_N-f\|_{L_2(\mathbb{R}_+)}=\|T_{\bar{g}_0,\bar{h}_0}(\varphi^N-\varphi)\|_{L_2(\mathbb{R}_+)}=\|\varphi^N-\varphi\|_{L_2 (\R_+)}$. Moreover, $\|\varphi^N-\varphi\|_{L_2 (\R_+)}\to 0$ when $N$ tends to $\infty$, which implies that $\|f_N-f\|_{L_2(\mathbb{R}_+)}\to 0$ means that a sequence of functions $\{f_N(x)\}$ converges to $f$ in the $L_2 (\R_+)$ norm. The proof of the second statement is completed.
\end{proof}
For real numbers $p, q>1$ , they are called conjugate indices (or Hölder conjugates) if
$
\frac{1}{p}+\frac{1}{q}=1.
$
Formally, we also define $q=\infty$ as conjugate to $p=1$ and vice versa.
Conjugate indices are used in Hölder's inequality, as well as Young's inequality for convolution products. If $p, q>1$ are conjugate indices, the spaces $L_p$ and $L_q$ are dual to each other. We will give the  condition for the boundedness of $(T_{g_0,h_0})$ operator on dual spaces through the following proposition
\begin{proposition}[Boundedness of $(T_{g_0,h_0})$ on dual spaces]\label{rem2} For case $\beta \in (0,1)$. Assume that
	$g_0\in L_1(\mathbb{R}_+)\cap L_2(\mathbb{R}_+)$ and  $h_0\in L_1^{0,\beta}(\mathbb{R}_+)\cap L_2^{0,\beta}(\mathbb{R}_+)$ simultaneously satisfy condition \eqref{eq6.2}, if $\left(1-d^2 / dx^2\right)\Phi(x,u,v,w)$  is a bounded function $\forall x\in \mathbb{R}_+$ then $(T_{g_0,h_0})$ defined by \eqref{eq6.1} is the bounded operator from $L_p(\mathbb{R}_+) \to L_q(\mathbb{R}_+)$, where $p\in[1,2]$ and $q$ is the conjugate exponent of $p$, that is,
	$\frac{1}{q}+\frac{1}{p}=1$.
\end{proposition}
\begin{proof}
	By Theorem \ref{thm6.2}, we infer that $(T_{g_0,h_0})$ is bounded from $L_2(\mathbb{R}_+)$ to $L_2(\mathbb{R}_+)$. From the above conditions, for any of the Lebesgue integrable functions over $\mathbb{R}_+$, we obtain
	$
	|\varphi(x)|=\left|\left(T_{g_0,h_0}f\right)(x)\right|<M
	$ is finite.
	This shows us that $(T_{g_0,h_0})$ is a bounded operator from  $L_1(\mathbb{R}_+)\to L_{\infty}(\mathbb{R}_+)$. By applying Riesz–Thorin's interpolation theorem \cite{Duoandikoetxea01}, we infer that the operator $(T_{g_0,h_0})$ is bounded from $L_p(\mathbb{R}_+)\to L_q(\mathbb{R}_+)$, where, $p$ and $q$ are conjugate exponents and $p$ is defined by $\frac1p=\frac{1-\alpha}{1}+\frac{\alpha}2=1-\frac{\alpha}2$ with $\alpha \in (0,1)$ and $1<p<2$. For the remaining cases, for instance, the case when $p=1$ which implies that $q=\infty$, and the case $p=2$, $q=2,$ have all been considered by us above. Therefore $(T_{g_0,h_0})$ is bounded from $L_p(\mathbb{R}_+)\to L_q(\mathbb{R}_+)$ when $p\in [1,2]$ with $\frac1p+\frac1q=1$.
\end{proof}	
\noindent By similar arguments, we also obtain this result for the operator $(T_{\bar{g}_0,\bar{h}_0})$ which is the inverse operator of $(T_{g_0,h_0})$ determined by \eqref{eq6.3}.

\section{Young and weighted $L_p$-norm inequalities for polyconvolution}\label{sec3}
\noindent The classical result for an upper bound estimation of the Fourier convolution
can be represented as follows:\\
$\lv f\underset{F}{*}g \lv_{L_s (\R)} \leq \lv f\lv_{L_p (\R)} \lv g\lv_{L_q (\R)},$
where $p,q,s > 1$ be real numbers such that
$\frac{1}{p} +\frac{1}{q}=1+\frac{1}{s}$ for any $f \in L_p (\R)$, $g\in L_q(\R)$. This result later became more widely known as Young's convolution inequality.
Subsequently, Adams and Fournier extended Young's inequality, transforming it into
a more general version of the Fourier convolution (see Theorem 2.24 in \cite{AdamsFournier2003sobolev}).
The above inequality can also be viewed as a consequence of Riesz's interpolation \cite{Duoandikoetxea01}. On the other hand, it is important to note that the Young inequality for Fourier convolution
is not valid in Hilbert space $L_2$. To address this limitation,
based on the general theory of reproducing kernels \cite{Saitoh2016},
Saitoh \cite{Saitoh2000} proposed a new inequality for the Fourier convolution in weighted
$L_{p}(\mathbb{R},|\rho_j|)$ Lebesgue spaces involving
non-vanishing functions $\rho_j$. We refer the reader to recent works \cite{TuanMMA2024,hoangtuan2017thaovkt,Tuan2023VKT,TuanHienPhuong}
for some other versions of these inequalities with applications to various transforms. Our purpose in this section is to derive a variant version of Young's inequality for operator \eqref{eq3.1}, establish the boundedness properties in a
three-parametric family of Lebesgues spaces as well as derive new weighted $L_p$-norm inequalities for this polyconvolution. The upper-bound constant in these inequalities is determined by using a special function, which is the Sech function. Sech is the hyperbolic secant function, which serves as the hyperbolic equivalent of the circular secant function used in trigonometry. It is defined as the reciprocal of the hyperbolic cosine function, expressed as $\mathrm{sech}\, t=\frac1{\cosh t}$. For real numbers, it is defined by considering twice the area between the axis and a ray through the origin that intersects the unit hyperbola. 
\subsection{The type of Young's inequalities and its direct corollary}
\begin{theorem}[Young-type theorem for operator \eqref{eq3.1}]\label{thm4.1}
Let $p,q,r$ and $s$ be real numbers in the interval $(1,\infty)$
such that $1/p + 1/q + 1/r+ 1/s=3$.	
For any functions $f\in L_p(\mathbb{R}_+)$, $g\in L_q(\mathbb{R}_+)$, $k\in L_s(\mathbb{R}_+)$ and $h\in L_r^{0,\beta}(\mathbb{R}_+)$ with $\beta \in (0,1]$, the following estimation holds	
	\begin{equation}\label{eq4.1}
	\left|\int_0^\infty \underset{F_s,F_c,K}{\ast}(f,g,h)(x)\cdot k(x) dx\right|\leq \mathcal{C}_1 \|f\|_{L_p(\mathbb{R}_+)} \|g\|_{L_q(\mathbb{R}_+)} \|h\|_{L_r^{0,\beta}(\mathbb{R}_+)},
	\end{equation}
	where $\mathcal{C}_1:=\sqrt{\frac{2}{\pi}}\left(\mathrm{sech}\, t\right)^{1-\frac1r}$ is the upper bound constant
	on the right-hand side of the inequality \eqref{eq4.1}.
\end{theorem}
\begin{proof} Without loss of generality, we assume that $p_1,q_1,r_1$ and $s_1$
	are the conjugate exponentials of $p,q,r$ and $s$, respectively. Since $1/p + 1/p_1=1$ (and similar identities for other constants),  together with the assumption of the theorem, then the following correlation between exponential numbers can be easily formulated
\begin{equation}\label{eq4.2}
\left\{\begin{array}{l}
\frac1{p_1}+\frac1{q_1}+\frac1{r_1}+\frac1{s_1}=1,\\
p\big(\frac{1}{q_1}+\frac{1}{r_1}+\frac{1}{s_1}\big)=q\big(\frac{1}{p_1}+\frac{1}{r_1}
+\frac{1}{s_1}\big)=r\big(\frac{1}{p_1}+\frac{1}{q_1}+\frac{1}{s_1}\big)
=s\big(\frac{1}{p_1}+\frac{1}{q_1}+\frac{1}{r_1}\big)=1,
\end{array}\right.
\end{equation}	
We set operators
	\begin{align*}
	E(x,u,v,w)&=|g(v)|^{\frac{q}{p_1}}|h(w)|^{\frac{r}{p_1}}|k(x)|^{\frac{s}{p_1}}|\Phi(x,u,v,w)|^{\frac{1}{p_1}}\in L_{p_1}(\mathbb{R}_+^4),\\
	F(x,u,v,w)&=|f(u)|^{\frac{p}{q_1}}|h(w)|^{\frac{r}{q_1}}|k(x)|^{\frac{s}{q_1}}|\Phi(x,u,v,w)|^{\frac{1}{q_1}}\in L_{q_1}(\mathbb{R}_+^4)	,\\
	H(x,u,v,w)&=|f(u)|^{\frac{p}{s_1}}|g(v)|^{\frac{q}{s_1}}|k(x)|^{\frac{r}{s_1}}|\Phi(x,u,v,w)|^{\frac{1}{s_1}}\in L_{s_1}(\mathbb{R}_+^4),\\
	G(x,u,v,w)&=|f(u)|^{\frac{p}{r_1}}|g(v)|^{\frac{q}{r_1}}|k(x)|^{\frac{s}{r_1}}|\Phi(x,u,v,w)|^{\frac{1}{r_1}}\in L_{r_1}(\mathbb{R}_+^4).
	\end{align*}
Due to condition \eqref{eq4.2}, we get $E(x,u,v,w) F(x,u,v,w) H(x,u,v,w) G(x,u,v,w)
=|f(u)||g(v)||h(w)||k(x)||\Phi(x,u,v,w)|$. Therefore   
$$\mathcal{I}:=\bigg|\displaystyle\int_0^\infty \underset{F_s,F_c,K}{\ast}(f,g,h)(x).k(x)dx\bigg|
\leq\frac{1}{2\sqrt{2\pi}}\int_{\mathbb{R}_+^4}|\Phi(x,u,v,w)||f(u)| |g(v)| |h(w)| |k(x)| du dv dw dx.$$ Hence $
\frac1{p_1}+\frac1{q_1}+\frac1{r_1}+\frac1{s_1}=1$,  by the four-function form of Hölder inequality, we deduce that
\begin{equation}\label{EFHG}
\begin{aligned}
\mathcal{I}\leq& \frac{1}{2\sqrt{2\pi}}\left(\int_{\mathbb{R}_+^4}|E(x,u,v,w)|^{p_1}dx du dv dw\right)^{\frac{1}{p_1}}\left(\int_{\mathbb{R}_+^4}|F(x,u,v,w)|^{q_1} dx du dv dw\right)^{\frac{1}{q_1}}\times\\
\times&\left(\int_{\mathbb{R}_+^4}|H(x,u,v,w)|^{s_1}dx du dv dw\right)^{\frac{1}{s_1}}\left(\int_{\mathbb{R}_+^4}|G(x,u,v,w)|^{r_1}dx du dv dw\right)^{\frac{1}{r_1}}\\
=&\frac{1}{2\sqrt{2\pi}}\|E\|_{L_{p_1}(\mathbb{R}_+^4)}\|F\|_{L_{q_1}(\mathbb{R}_+^4)}\|H\|_{L_{s_1}(\mathbb{R}_+^4)}\|G\|_{L_{r_1}(\mathbb{R}_+^4)}.
\end{aligned}
\end{equation}
Based on the assumption of $f\in L_p(\mathbb{R}_+)$, $g\in L_q(\mathbb{R}_+)$, $k\in L_s(\mathbb{R}_+)$ and $h\in L_r^{0,\beta}(\mathbb{R}_+)$, by using Fubini's theorem and combining \eqref{eq3.3},  we obtain  $L_{p_1} (\R^4_+)$-norm estimation for the operator $E$ as follows:
	
	$$\begin{aligned}
\|E\|_{L_{p_1}(\mathbb{R}_+^4)}^{p_1}&=\int_{\mathbb{R}_+^4}|g(v)|^q|h(w)|^r|k(x)|^s|\Phi(x,u,v,w)| dv dw dx du\\
&=\int_0^\infty |g(v)|^q\left\{\int_0^\infty|h(w)|^r\bigg[\int_0^\infty|k(x)|^s\bigg(\int_0^\infty|\Phi(x,u,v,w)|du\bigg)dx\bigg]dw\right\}dv\\
&\leq \int_0^\infty |g(v)|^q\left\{\int_0^\infty 4K_0(w)|h(w)|^r\left(\int_0^\infty |k(x)|^sdx\right)dw\right\}dv\\
&\leq 4\left(\int_0^\infty|g(v)|^q dv\right)\left(\int_0^\infty K_0(\beta w)|h(w)|^r dw\right)\left(\int_0^\infty |k(x)|^s dx\right) \\
&=4\|g(v)\|_{L_q(\mathbb{R}_+)}^q \|h\|^r_{L_r^{0,\beta}(\mathbb{R}_+)}\|k\|_{L_s(\mathbb{R}_+)}^s\ \text{with}\ \beta \in (0,1].
\end{aligned}$$
This yields an $L_{p_1}(\mathbb{R}_+^4)$-norm estimation for the operator $E$ as
\begin{equation}\label{eq4.4}
\|E\|_{L_{p_1}(\mathbb{R}_+^4)}\leq 4^{\frac{1}{p_1}}\|g\|_{L_q(\mathbb{R}_+)}^{\frac{q}{p_1}}\|h\|_{L_r^{0,\beta}(\mathbb{R}_+)}^{\frac{r}{p_1}}\|k\|_{L_s(\mathbb{R}_+)}^{\frac{s}{p_1}}.
\end{equation}
By similar arguments used to derive $L_{p_1}(\mathbb{R}_+^4)$-norm estimation for $E$, we also have
	\begin{equation}\label{eq4.5}
	\|F\|_{L_{q_1}(\mathbb{R}_+^4)}\leq 4^{\frac{1}{q_1}}\|f\|_{L_p(\mathbb{R}_+)}^{\frac{p}{q_1}}\|h\|_{L_r^{0,\beta}(\mathbb{R}_+)}^{\frac{r}{q_1}}\|k\|_{L_s(\mathbb{R}_+)}^{\frac{s}{q_1}},\ \textup{and}\ \ \|H\|_{L_{s_1}(\mathbb{R}_+^4)}\leq 4^{\frac{1}{s_1}}\|f\|_{L_p(\mathbb{R}_+)}^{\frac{p}{s_1}}\|g\|_{L_q(\mathbb{R}_+)}^{\frac{q}{s_1}}\|h\|_{L_r^{0,\beta}(\mathbb{R}_+)}^{\frac{r}{s_1}}.
	\end{equation}
	To give an estimate for $G$,  it is easy to first see that $\int_0^\infty|\Phi(x,u,v,w)|dw \leq \frac{4}{\cosh t}$ due to \eqref{eq3.3}. This means that
	\begin{equation}\label{eq4.7}
	\|G\|_{L_{r_1}(\mathbb{R}_+^4)}\leq \left(\frac{4}{\cosh t}\right)^{\frac{1}{r_1}}\|f\|_{L_p(\mathbb{R}_+)}^{\frac{p}{r_1}}\|g\|_{L_q(\mathbb{R}_+)}^{\frac{q}{r_1}}\|k\|_{L_s(\mathbb{R}_+)}^{\frac{s}{r_1}}.
	\end{equation}
	Coupling  \eqref{eq4.4},\eqref{eq4.5},\eqref{eq4.7},  and by taking condition \eqref{eq4.2} into
	account, we readily obtain 
	\begin{align}\label{eq4.8}
	\|E\|_{L_{p_1}(\mathbb{R}_+^4)}\|F\|_{L_{q_1}(\mathbb{R}_+^4)}\|H\|_{L_{s_1}(\mathbb{R}_+^4)}\|G\|_{L_{r_1}(\mathbb{R}_+^4)}
	\leq 4\left(\frac{1}{\cosh t}\right)^{1-\frac{1}{r}}\|f\|_{L_p(\mathbb{R}_+)}\|g\|_{L_q(\mathbb{R}_+)}\|h\|_{L_r^{0,\beta}(\mathbb{R}_+)}\|k\|_{L_s(\mathbb{R}_+)}.
	\end{align}
	 Finally, with  a combination of \eqref{EFHG} and \eqref{eq4.8}, we obtain an estimation as in the conclusion of this theorem with upper bound constant $\mathcal{C}_1= \frac{1}{2\sqrt{2\pi}}. 4\left(\frac{1}{\cosh t}\right)^{1-\frac{1}{r}}= \sqrt{\frac{2}{\pi}}\left(\mathrm{sech}\, t\right)^{1-\frac1r}$ with $\mathrm{sech}\, t=\frac1{\cosh t}$.
\end{proof}
\noindent In a special case, for $k(x)$ defined by \eqref{eq3.1}, the following Young-type inequality is a direct consequence of Theorem \ref{thm4.1}.
\begin{corollary}\label{cor4.1}{\rm(Young-type inequality)}
	 Suppose that $p,q,r,s\in[1,\infty)$ satisfying $\frac{1}{p}+\frac{1}{q}+\frac{1}{r}=2+\frac{1}{s}$. Let $f\in L_p(\mathbb{R}_+), g\in L_q(\mathbb{R}_+), h\in L_r^{0,\beta}(\mathbb{R}_+)$ with $0<\beta\leq 1$, then and the polyconvolution \eqref{eq3.1} is well-defined and belongs to $L_s(\mathbb{R}_+)$. Hence, the following estimate always holds  true:
	\begin{equation}\label{eq4.10}
	\big\|\underset{F_s,F_c,K}{\ast}(f,g,h)\big\|_{L_s(\mathbb{R}_+)}\leq \mathcal{C}_1 \|f\|_{L_p(\mathbb{R}_+)}\|g\|_{L_q(\mathbb{R}_+)}\|h\|_{L_r^{0,\beta}(\mathbb{R}_+)}. 
	\end{equation} 
\end{corollary}
\begin{proof}
	For case $p=q=r=s=1$, this result has been proved in Theorem ~\ref{thm3.1}. Thus, we only need to confirm this theorem holds for cases when $p,q,r,s>1.$ Indeed, 
assume that $s_1$ is the conjugate exponent of $s$, i.e. $1/s +1/s_1 =1$. From the assumptions of
Corollary \ref{cor4.1}, we have $\frac1p+\frac1q+\frac1r+\frac1{s_1}=3$, which shows the numbers $p$, $q$, $r$, and $s_1$ satisfy the condition of Theorem \ref{thm4.1} (with the role of $s$ being replaced by  $s_1$). Therefore, with $f\in L_p(\mathbb{R}_+)$, $g\in L_q(\mathbb{R}_+)$ and $h\in L_r^{0,\beta}(\mathbb{R}_+)$, the linear operator 
	$
	\mathscr{T}_k:=\int_0^\infty \underset{F_s,F_c,K}{\ast}(f,g,h)(x). k(x) dx
	$
	is bounded $L_{s_1}(\mathbb{R}_+)$. Consequently, by using the Riesz representation theorem \cite{Stein1971Weiss} then   $\underset{F_s,F_c,K}{\ast}(f,g,h)$ belongs to $L_{s_1}(\mathbb{R}_+)$. To continue the proof, we choose
	$
	k(x):=\mathrm{sign}\big[\underset{F_s,F_c,K}{\ast}(f,g,h)(x)\big]^s \times\big[\underset{F_s,F_c,K}{\ast}(f,g,h)(x)\big]^{\frac{s}{s_1}}.
	$
	Therefore, $k\in L_{s_1}(\mathbb{R}_+)$ and obtain the norm
	$\|k\|_{L_{s_1}(\mathbb{R}_+)}=\big\|\underset{F_s,F_c,K}{\ast}f,g,h\big\|_{L_s(\mathbb{R}_+)}^{s/s_1}.$ By 
	applying the inequality \eqref{eq4.1}  to such function $k(x)$ we deduce that
	$$\begin{aligned}
	\big\|\underset{F_s,F_c,K}{\ast}f,g,h\big\|_{L_s(\mathbb{R}_+)}^s&=\left|\int_0^\infty \underset{F_s,F_c,K}{\ast}(f,g,h)(x)dx\right|^s = \left| \int_0^\infty \underset{F_s,F_c,K}{\ast}(f,g,h)(x)\cdot k(x) dx\right|\\
	&\leq \sqrt{\frac{2}{\pi}}\left(\mathrm{sech}\, t\right)^{1-\frac1r} \|f\|_{L_p(\mathbb{R}_+)} \|g\|_{L_q(\mathbb{R}_+)}\|h\|_{L_r^{0,\beta}(\mathbb{R})}\|k\|_{L_{s_1}(\mathbb{R}_+)} \\
	&=\sqrt{\frac{2}{\pi}}\left(\mathrm{sech}\, t\right)^{1-\frac1r} \|f\|_{L_p(\mathbb{R}_+)} \|g\|_{L_q(\mathbb{R}_+)}\|h\|_{L_r^{0,\beta}(\mathbb{R})} \big\|\underset{F_s,F_c,K}{\ast}f,g,h\big\|_{L_s(\mathbb{R}_+)}^{s/s_1}, 
	\end{aligned}$$
or equivalent $\forall \beta \in (0,1]$, then $\big\|\underset{F_s,F_c,K}{\ast}f,g,h\big\|_{L_s(\mathbb{R}_+)}^{s-\frac{s}{s_1}}\leq \sqrt{\frac{2}{\pi}}\left(\mathrm{sech}\, t\right)^{1-\frac1r} \|f\|_{L_p(\mathbb{R}_+)}\|g\|_{L_q(\mathbb{R}_+)}\|h\|_{L_r^{0,\beta}(\mathbb{R}_+)}.$
	Since $\frac1{s}+\frac1{s_1}=1$, which implies that $s-\frac{s}{s_1}=1$ and from there we arrive at the conclusion of the corollary.
\end{proof}
From the condition of Corollary ~\ref{cor4.1} and using the estimate \eqref{eq4.10} we get 
\begin{remark}\label{rmk2}
\textup{	If $f,g$ are $L_1$-Lebesgue integrable functions over $\R_+$ and $h\in L_s^{0,\beta}(\mathbb{R}_+)$, then the estimate holds true}
	\begin{equation}\label{eq4.12a}
	\big\| \underset{F_s,F_c,K}{\ast} (f,g,h)\big\|_{L_s(\mathbb{R}_+)} \leq \mathcal{C}_1 \|f\|_{L_1(\mathbb{R}_+)} \|g\|_{L_1(\mathbb{R}_+)}\|h\|_{L_s^{0,\beta}(\mathbb{R}_+)},\forall \beta \in (0,1], \textup{and}\   s\geq 1.
	\end{equation}
\end{remark}
\noindent What can be concluded for the case $s=\infty$ is a natural and essential question. The boundedness of the operator \eqref{eq3.1} in  $L_{\infty}(\R_+)$ is presented in the following proposition.
\begin{proposition}\label{thm4.2}
	Let $p,q,r\in (1,\infty)$ such that $\frac1p+\frac1q+\frac1r=2$. For any $f\in L_p(\mathbb{R}_+), g\in L_q(\mathbb{R}_+), h\in L_r^{0,\beta}(\mathbb{R}_+)$ $\forall \beta \in (0,1]$ then the operator \eqref{eq3.1} is a bounded function $\forall x\in \mathbb{R}_+$ and there is the following estimate
	\begin{equation}
	\label{eq4.13} \big\|\underset{F_s,F_c,K}{\ast}f,g,h\big\|_{L_\infty(\mathbb{R}_+)}\leq \mathcal{C}_1 \|f\|_{L_p(\mathbb{R}_+)}\|g\|_{L_q(\mathbb{R}_+)}\|h\|_{L_r^{0,\beta}(\mathbb{R}_+)},
	\end{equation}
	where the norm on $L_\infty(\mathbb{R}_+)$ is defined via essential supremum \textquotedblleft ess $\sup $\textquotedblright i.e. $\|\varphi\|_{L_\infty(\mathbb{R}_+)}=ess\sup\limits_{x\in \mathbb{R}_+} |\varphi(x)|$. 
\end{proposition}
\begin{proof}
	Suppose that $p_1$, $q_1$, $r_1$ be the conjugate exponentials of $p$, $q$, $r$  respectively, together with the assumption $\frac1p+\frac1q+\frac1r=2$ of Proposition \eqref{thm4.2}, we deduce that $\frac{1}{p_1}+\frac{1}{q_1}+\frac{1}{r_1}=1$ and $\frac{p}{q_1}+\frac{p}{r_1}=\frac{q}{q_1}+\frac{q}{r_1}=\frac{r}{p_1}+\frac{r}{q_1}=1.$
	We set
	\begin{align*}
	N(u,v,w)&=|g(v)|^{\frac{q}{p_1}} |h(w)|^{\frac{r}{p_1}} |\Phi(x,u,v,w)|^{\frac{1}{p_1}}\in L_{p_1}(\mathbb{R}_+^3),\\
	T(u,v,w)&=|f(u)|^{\frac{p}{q_1}} |h(w)|^{\frac{r}{q_1}} |\Phi(x,u,v,w)|^{\frac{1}{q_1}}\in L_{q_1}(\mathbb{R}_+^3),\\
	M(u,v,w)&=|f(u)|^{\frac{p}{r_1}} |g(v)|^{\frac{q}{r_1}} |\Phi(x,u,v,w)|^{\frac{1}{r_1}}\in L_{r_1}(\mathbb{R}_+^3).
	\end{align*} This yields $N(u,v,w).T(u,v,w). M(u,v,w) = |f(u)| |g(v)| |h(w)| |\Phi(x,u,v,w)| $ by the above setup because of the conjugate exponent constraint. Therefore
	\begin{align}
	\big|\underset{F_s,F_c,K}{\ast}(f,g,h)(x)\big|&\leq \frac{1}{2\sqrt{2\pi}}\int_{\mathbb{R}_+^3}|\Phi(x,u,v,w)| |f(u)| |g(v)| |h(w)| du dv dw \nonumber\\
	&=\frac{1}{2\sqrt{2\pi}}\int_{\mathbb{R}_+^3} N(u,v,w)T(u,v,w)M(u,v,w) du dv dw.\label{eq4.16}
	\end{align}
	However, $\frac{1}{p_1}+\frac{1}{q_1}+\frac{1}{r_1}=1$. Using the H\"older's inequality into \eqref{eq4.16} we obtain
	\begin{equation}\label{eq4.17}
	\begin{aligned}
	&\big|\underset{F_s,F_c,K}{\ast}(f,g,h)(x)\big|\\&\leq \frac{1}{2\sqrt{2\pi}}\left\{\int_{\mathbb{R}_+^3}|N(u,v,w)|^{p_1}du dv dw\right\}^{\frac{1}{p_1}}
	\times \left\{\int_{\mathbb{R}_+^3}|T(u,v,w)|^{q_1}dudvdw\right\}^{\frac{1}{q_1}}\times \left\{\int_{\mathbb{R}_+^3}|M(u,v,w)^{r_1}du dv dw\right\}^{\frac1{r_1}}\\
	&=\frac{1}{2\sqrt{2\pi}}\|N\|_{L_{p_1}(\mathbb{R}_+^3)} \|T\|_{L_{q_1}(\mathbb{R}_+^3)} \|M\|_{L_{r_1}(\mathbb{R}_+^3)}
	\end{aligned}
	\end{equation}
	Similar to the proof of Theorem \ref{thm4.1} we infer the norm estimates of $N$, $T$, $M$ and have form as:
	\begin{equation}\label{eq4.18}
	\|N\|_{L_{p_1}(\mathbb{R}_+^3)}\leq 4^{\frac{1}{p_1}}\|g\|_{L_q(\mathbb{R}_+)}^{\frac{q}{p_1}}\|h\|_{L_r^{0,\beta}(\mathbb{R}_+)}^{\frac{r}{p_1}},\ \textup{and}\ \  \|T\|_{L_{q_1}(\mathbb{R}_+^3)}\leq 4^{\frac{1}{q_1}}\|f\|_{L_p(\mathbb{R}_+)}^{\frac{p}{q_1}}\|h\|_{L_r^{0,\beta}(\mathbb{R}_+)}^{\frac{r}{q_1}},
	\end{equation}
	and
	\begin{equation}\label{eq4.20}
	\|M\|_{L_{r_1}(\mathbb{R}_+^3)}\leq \left(\frac4{\cosh t}\right)^{\frac{1}{r_1}}\|f\|_{L_p(\mathbb{R}_+)}^{\frac{p}{r_1}}\|g\|_{L_q(\mathbb{R}_+)}^{\frac{q}{r_1}}.
	\end{equation}
	Coupling \eqref{eq4.18},\eqref{eq4.20} and $\forall f\in L_p(\mathbb{R}_+)$, $g\in L_q(\mathbb{R}_+)$, $h\in L_r^{0,\beta}(\mathbb{R}_+)$ with $\forall \beta \in (0,1]$, we obtain
	\begin{equation}\label{eq4.21}
	\begin{aligned}
	\|N\|_{L_{p_1}(\mathbb{R}_+^3)}\|T\|_{L_{q_1}(\mathbb{R}_+^3)}\|M\|_{L_{r_1}(\mathbb{R}_+^3)}
	&\leq  4^{\left(\frac1{p_1}+\frac1{q_1}+\frac1{r_1}\right)}\left(\frac1{\cosh t}\right)^{\frac1{r_1}} \|f\|_{L_p(\mathbb{R}_+)}^{\left(\frac{p}{q_1}+\frac{p}{r_1}\right)} \|g\|_{L_q(\mathbb{R}_+)}^{\left(\frac{q}{p_1}+\frac{q}{r_1}\right)} \|h\|_{L_r^{0,\beta}(\mathbb{R}_+)}^{\left(\frac{r}{p_1}+\frac{r}{q_1}\right)}\\
	&=4\left(\mathrm{sech}\,t\right)^{1-\frac1r}\|f\|_{L_p(\mathbb{R}_+)}\|g\|_{L_q(\mathbb{R}_+)}\|h\|_{L_r^{0,\beta}(\mathbb{R}_+)}<\infty.
	\end{aligned}\end{equation}
	Finally, the combination of  \eqref{eq4.17} and \eqref{eq4.21}  gives  \eqref{eq3.1} as a bounded operator in $L_{\infty}$ and $ess\sup\limits_{x\in\mathbb{R}_+} \big|\underset{F_s,F_c,K}{\ast}(f,g,h)(x)\big|$ which is finite. The proof is completed.
\end{proof}
\noindent Theorem \ref{thm3.1} and inequalities \eqref{eq4.10},\eqref{eq4.13} show that \eqref{eq3.1}  is valid and well-defined with $s \in [1,\infty]$. Moreover, for any $p,q,r\geq1$, these have also shown that \eqref{eq3.1} is a bounded multilinear map from $L_p (\R_+)\times L_q (\R_+)\times L^{0,\beta}_r (\R_+)$ to $L_s (\R_+)$ such that $(f,g,h)\mapsto \underset{F_s,F_c,K}{\ast}(f,g,h)$.
\subsection{Boundedness on index spaces $L_s^{\alpha_1,\beta_1,\gamma_1}(\mathbb{R}_+)$}\label{subsec32}
\begin{theorem}\label{thm4.3}
Suppose that $p,q,r$ be real numbers in $(1,\infty)$ such that $\frac1p+\frac1q+\frac1r=2$. For any function $f\in L_p(\mathbb{R}_+)$, $g\in L_q(\mathbb{R}_+)$, $h\in L_r^{0,\beta}(\mathbb{R}_+)$ with  $\beta \in (0,1]$, then operator \eqref{eq3.1} is well-defined as the continuous function and belongs to $L_s^{\alpha_1,\beta_1,\gamma_1}(\mathbb{R}_+)$ with  $s\geq 1$ and three-parametric $\alpha_1>-1$, $\beta_1,\gamma_1>0$. Furthermore
\begin{equation}\label{eq4.23}
\big\|\underset{F_s,F_c,K}{\ast} f,g,h\big\|_{L_s^{\alpha_1,\beta_1,\gamma_1}(\mathbb{R}_+)} \leq \mathcal{C}_2\|f\|_{L_p(\mathbb{R}_+)} \|g\|_{L_q(\mathbb{R}_+)}\|h\|_{L_r^{0,\beta}(\mathbb{R}_+)},
\end{equation}
where $\mathcal{C}_2:=\sqrt{\frac{2}{\pi}}(\mathrm{sech}\,t)^{1-\frac1r} \left(\frac1{\gamma_1}\right)^{\frac1s} \beta_1^{\frac{-\alpha_1+1}{\gamma_1s}}\Gamma^{\frac1s}\left(\frac{\alpha_1+1}{\gamma_1}\right)$ is the upper bound coefficient on the right-hand side of inequality \eqref{eq4.23}. Here $
L_s^{\alpha_1,\beta_1,\gamma_1} (\mathbb{R}_+):=\left\{\varphi(x)\ \textup{defined in}\ \mathbb{R}_+ : \int_0^\infty |\varphi(x)|^s x^{\alpha_1} e^{-\beta_1 x^{\gamma_1}}dx<\infty\right\}
$ is the three-parametric family of Lebesgue spaces defined similarly as in \cite{Yakubovich1996index} 
with the norm
$
\|\varphi\|_{L_s^{\alpha_1,\beta_1,\gamma_1}(\mathbb{R}_+)}=\left\{\int_0^\infty |\varphi(x)|^s x^{\alpha_1} e^{-\beta_1 x^{\gamma_1}}dx\right\}^\frac{1}{s}.
$
\end{theorem}
\begin{proof}
	From \eqref{eq4.13}, we deduce
	$\big|\underset{F_s,F_c,K}{\ast}(f,g,h)(x)\big| \leq \sqrt{\frac{2}{\pi}}(\mathrm{sech}\, t)^{1-\frac1r} \|f\|_{L_p(\mathbb{R}_+)} \|g\|_{L_1(\mathbb{R}_+)} \|h\|_{L_r^{0,\beta}(\mathbb{R}_+)}
	\leq \mathcal{M}$ is finite, where $\mathcal{M}$ is a positive constant. According to the formula 3.381.4 in \cite{gradshteyn2014Ryzhik} with $\alpha_1>-1$, $\beta_1,\gamma_1>0$, we infer that
	$
	\int_0^\infty x^{\alpha_1} e^{-\beta_1x^{\gamma_1}} dx = \frac1{\gamma_1}\beta_1^{\frac{-\alpha_1+1}{\gamma_1}}\Gamma\left(\frac{\alpha_1+1}{\gamma_1}\right).
	$ Therefore
	$$\int_0^\infty \big|\underset{F_c,F_s,K}{\ast}(f,g,h)(x)\big|^s x^{\alpha_1} e^{-\beta_1 x^{\gamma_1}} dx\leq \mathcal{M}^s\int_0^\infty x^{\alpha_1} e^{-\beta_1 x^{\gamma_1}}dx
	=\mathcal{M}^s~ \frac{1}{\gamma_1}\beta_1^{\frac{-\alpha_1+1}{\gamma_1}}\Gamma\left(\frac{\alpha_1+1}{\gamma_1}\right) <\infty,$$ which implies that polyconvolution \eqref{eq3.1} belongs to $L_s^{\alpha_1,\beta_1,\gamma_1}(\mathbb{R}_+)$ and we obtain the below estimate
	$$\begin{aligned}
	\big\|\underset{F_s,F_c,K}{\ast} f,g,h\big\|_{L_s^{\alpha_1,\beta_1,\gamma_1}(\mathbb{R}_+)} \leq \sqrt{\frac{2}{\pi}}(\mathrm{sech}\,t)^{1-\frac1r}~ \left(\frac1{\gamma_1}\right)^{\frac1s}\beta_1^{\frac{-\alpha_1+1}{\gamma_1s}}\Gamma^{\frac1s}\left(\frac{\alpha_1+1}{\gamma_1}\right)
\|f\|_{L_p(\mathbb{R}_+)} \|g\|_{L_q(\mathbb{R}_+)} \|h\|_{L_r^{0,\beta}(\mathbb{R}_+)}.
	\end{aligned}$$
	This leads to  $\mathcal{C}_2=\sqrt{\frac{2}{\pi}}(\mathrm{sech}\,t)^{1-\frac1r} \left(\frac1{\gamma_1}\right)^{\frac1s} \beta_1^{\frac{-\alpha_1+1}{\gamma_1s}}\Gamma^{\frac1s}\left(\frac{\alpha_1+1}{\gamma_1}\right)=\mathcal{C}_1 \bigg[\left(\frac1{\gamma_1}\right)^{\frac1s} \beta_1^{\frac{-\alpha_1+1}{\gamma_1s}}\Gamma^{\frac1s}\left(\frac{\alpha_1+1}{\gamma_1}\right)\bigg] $, where $\Gamma(x)$ is an Euler's gamma-function \cite{Yaku94Luchko}.
\end{proof}

\subsection{The weighted inequality in classes $L_p(\mathbb{R}_+,\rho)$ space}
\begin{theorem}[Saitoh-type inequality for operator \eqref{eq3.1}]\label{thm5.1} Assume
	that $\rho_j$ are non-vanishing positive functions for all $x\in \R_+$,  then 
	 $\underset{F_s,F_c,K}{\ast}(\rho_1,\rho_2,\rho_3)(x)$ given by the form of \eqref{eq3.1}  is well-defined. Moreover $\forall F_j\in L_p(\mathbb{R}_+,\rho_j),$ $j=\overline{1,3}$ with $p>1$  the following weighted $L_p (\R_+)$-norm inequality holds
	\begin{equation}\label{eq5.1}
	\big\|\underset{F_s,F_c,K}{\ast}(F_1\rho_1,F_2\rho_2,F_3\rho_3)\times \underset{F_s,F_c,K}{\ast}(\rho_1,\rho_2,\rho_3)^{\frac{1}{p}-1}\big\|_{L_p(\mathbb{R}_+)}
	\leq \left(\sqrt{\frac{2}{\pi}}K_0(w)\right)^{\frac1p} \prod_{i=j}^3\|F_j\|_{L_p(\mathbb{R}_+,\rho_j)}.
	\end{equation}
	Here the $K_0(\omega)$ function is calculated via  $K_0(\omega)=\int_0^\infty e^{-\omega\cosh u}du, \forall w>0$
\end{theorem}
\begin{proof}
	According to \eqref{eq3.1}, we have 
	\begin{equation}\label{eq5.2}
\begin{aligned}
&\big\|\underset{F_s,F_c,K}{\ast}(F_1\rho_1,F_2\rho_2,F_3\rho_3)\times \underset{F_s,F_c,K}{\ast}(\rho_1,\rho_2,\rho_2)^{\frac{1}{p}-1}\big\|_{L_p(\mathbb{R}_+)}^p\\
&=\int_0^\infty \left|\underset{F_s,F_c,K}{\ast}(F_1\rho_1,F_2\rho_2,F_3\rho_3)(x)\right|^p \left|\underset{F_s,F_c,K}{\ast}(\rho_1,\rho_2,\rho_3)(x)\right|^{1-p}dx\\
&\leq \left(\frac{1}{2\sqrt{2\pi}}\right)^p\left(\frac{1}{2\sqrt{2\pi}}\right)^{p-1}\int_0^\infty \bigg\{ \bigg(\int_{\mathbb{R}_+^3} |\Phi(x,u,v,w)| |(F_1\rho_1)(u)| |(F_2\rho_2)(v)||(F_3\rho_3)|(w)| du dv dw\bigg)^p\\&\times \left(\int_{\mathbb{R}_+^3}|\Phi(x,u,v,w)| \rho_1(u) \rho_2(v)\rho_3(w) du dv dw\right)^{1-p}\bigg\}dx,
\end{aligned}	
	\end{equation}
	where kernel $\Phi(x,u,v,w)$ is defined by \eqref{eq3.2}
	 Applying the Hölder inequality for the conjugate pair  $p$, $q$, we have
	\begin{equation}\label{eq5.3}
	\begin{aligned}
	&\int_{\mathbb{R}_+^3} |\Phi(x,u,v,w)| |(F_1\rho_1)(u)| |(F_2\rho_2)(v)| |(F_3\rho_3)(w)| du dv dw\\
	&\leq \left(\int_{\mathbb{R}^3} |\Phi(x,u,v,w)| |F_1(u)|^p \rho_1(u)|F_2(v)|^p\rho_2(v) |F_3(w)|^p\rho_3(w) du dv dw\right)^{\frac1p}\times\\&\times\left(\int_{\mathbb{R}_+^3} |\Phi(x,u,v,w)|\rho_1(u)\rho_2(v)\rho_3(w) du dv dw\right)^{\frac1q}.
	\end{aligned}
	\end{equation}
	It can be deduced from  \eqref{eq5.2} and \eqref{eq5.3} that
	\begin{equation}\label{eq5.4}
	\begin{aligned}
	&\big\|\underset{F_s,F_c,K}{\ast}(F_1\rho_1,F_2\rho_2,F_3\rho_3)\times \underset{F_s,F_c,K}{\ast}(\rho_1,\rho_2,\rho_3)^{\frac1p-1}\big\|_{L_p(\mathbb{R}_+)}^{\rho}\\
	&\leq \frac{1}{2\sqrt{2\pi}}\int_0^\infty\bigg\{\left(\int_0^\infty |\Phi(x,u,v,w)| |F_1(u)|^p\rho_1(u)|F_2(v)|^p\rho_2(v)|F_3(w)|^p\rho_3(w) du dv dw\right)\times\\
	&\times \left(\int_{\mathbb{R}_+^3}|\Phi(x,u,v,w)| \rho_1(u)\rho_2(v)\rho_3(w) du dv dw\right)^{\frac{p}{q}+1-p}\bigg\}dx.
	\end{aligned}
	\end{equation}
	Since $p$ and $q$ are a conjugate pair i.e. $\frac1p+\frac1q=1$, we have $\frac{p}{q}+1-p=0$.  By
	the assumption $F_j \in L_p (\R_+, \rho_j)$ using Fubini's Theorem for the integral \eqref{eq5.4} combining \eqref{eq3.3}, we infer
	$$\begin{aligned}
	&\left\|\underset{F_s,F_c,K}{\ast}(F_1\rho_1,F_2\rho_2,F_3\rho_3)\times \underset{F_s,F_c,K}{\ast}(\rho_1,\rho_2,\rho_3)^{\frac1p-1}\right\|_{L_p(\mathbb{R}_+)}^{\rho}\\
	&\leq \frac{1}{2\sqrt{2\pi}} \int_{\mathbb{R}_+^4}|\Phi(x,u,v,w)| |F_1(u)|^p \rho_1(u) |F_2(v)|^p\rho_2(v)|F_3(w)|^p\rho_3(w) du dv dw dx\nonumber\\
	&=\frac{1}{2\sqrt{2\pi}}\left(\int_0^\infty |\Phi(x,u,v,w)|dx\right)\left(\int_0^\infty|F_1(u)|^p\rho_1(u)du\right)\left(\int_0^\infty|F_2(v)|^p\rho_2(v)dv\right)\left(\int_0^\infty |F_3(w)|^p\rho_3(w)dw\right)\\
	&\leq \sqrt{\frac{2}{\pi}}K_0(w)\|F_1\|_{L_p(\mathbb{R}_+,\rho_1)}^p\|F_2\|_{L_p(\mathbb{R}_+,\rho_2)}^p\|F_3\|_{L_p(\mathbb{R}_+,\rho_3)}^p.
	\end{aligned}$$
\end{proof}

\noindent Next, we consider several special cases regarding weighted functions $\rho_j\in L_p(\mathbb{R}_+,\rho_i)$ space. If one of the weight functions $\rho_1(x)$ or
$\rho_2(x)$ or $\rho_3(x)$  is equal to $1$, for instance;\\ \textbf{Case 1:} Let $\rho_1(x)\equiv 1$ and $\rho_2 (x), \rho_3(x)$ be non-negative functions that belong to $L_1(\mathbb{R})$,
by combining  \eqref{eq3.1}, \eqref{eq3.3} with Fubini's theorem we can conclude that
\begin{align*}
&\big|\underset{F_s,F_c,K}{\ast}(1,\rho_2,\rho_3)(x)\big|\leq \frac{1}{2\sqrt{2\pi}} \int_{\mathbb{R}_+^3}|\Phi(x,u,v,w)|\rho_2(v)\rho_3(w) du dv dw\\
&\leq \frac{1}{2\sqrt{2\pi}}\left(\int_0^\infty|\Phi(x,u,v,w)|du\right)\left(\int_0^\infty \rho_2(v)dv\right)\left(\int_0^\infty \rho_3(w)dw\right)	
=\sqrt{\frac{2}{\pi}}K_0(w)\|\rho_2\|_{L_1(\mathbb{R}_+)}\|\rho_3\|_{L_1(\mathbb{R}_+)}
\end{align*}
which is finite.
This means that  $\underset{F_s,F_c,K}{\ast}(1,\rho_1,\rho_2)$ is well-defined and is a bounded
function  $\forall x\in\mathbb{R}_+$. Therefore, we infer that estimates
$
\big|\underset{F_s,F_c,K}{\ast}(1,\rho_1,\rho_2)^{1-\frac1p}(x)\big|\leq \left(\sqrt{\frac{2}{\pi}}K_0(w)\right)^{1-\frac1p} \|\rho_2\|_{L_1(\mathbb{R}_+)}^{1-\frac1p}\|\rho_3\|_{L_1(\mathbb{R}_+)}^{1-\frac1p}.
$
Combining with the inequality \eqref{eq5.1}, we obtain the following corollary.
\begin{corollary}\label{cor5.1}
	Let $\rho_1(x) \equiv 1\ \forall x \in \R_+$, $0<\rho_2,\rho_3\in L_1(\mathbb{R}_+)$. If  $F_1 \in L_p(\mathbb{R}_+)$ and $F_j\in L_p(\mathbb{R}_+,\rho_j),j=2,3$, then 
	\begin{equation}\label{eq5.6}
	\big\|\underset{F_s,F_c,K}{\ast}(F_1,F_2\rho_2,F_3\rho_3)\big\|_{L_p(\mathbb{R}_+)} \leq \mathcal{C}_3\|F_1\|_{L_p(\mathbb{R}_+)}\|F_2\|_{L_p(\mathbb{R}_+,\rho_2)}\|F_3\|_{L_p(\mathbb{R}_+,\rho_3)},
	\end{equation}
	 where the coefficient
	$
	\mathcal{C}_3:=\sqrt{\frac{2}{\pi}}K_0(w)\|\rho_2\|_{L_1(\mathbb{R}_+)}^{1-\frac1p}\|\rho_3\|_{L_1(\mathbb{R}_+)}^{1-\frac1p},
	$ with $p>1.$
\end{corollary}
\noindent We illustrate the estimation \eqref{eq5.6} by choosing weighted functions $\rho_1=1$, $\rho_2=e^{-t}$, $\rho_3=e^{-2t}\in L_1(\mathbb{R}_+)$. Hence $\big|\underset{F_s,F_c,K}{\ast}\left(1,e^{-t},e^{-2t}\right)(x)\big|\leq \sqrt{\frac{2}{\pi}}K_0(w)$ and $\big|\underset{F_s,F_c,K}{\ast}\left(1,e^{-t},e^{-2t}\right)^{1-\frac1p}(x)\big|\leq \left(\sqrt{\frac{2}{\pi}}K_0(w)\right)^{1-\frac1p}, w>0.$
Applying \eqref{eq5.1} we get
$
\big\|\underset{F_s,F_c,K}{\ast}\left(F_1,F_2e^{-t},F_2e^{-2t}\right)\big\|_{L_p(\mathbb{R}_+)}\leq \sqrt{\frac{2}{\pi}}K_0(w)\|F_1\|_{L_p(\mathbb{R}_+)} \|F_2\|_{L_p(\mathbb{R}_+,e^{-t})} \|F_3\|_{L_p(\mathbb{R}_+,e^{-2t})}.
$
For case $\rho_2 (x)\equiv 1$, $0<\rho_1,\rho_3\in L_1(\mathbb{R}_+)$, we also get a similar inequality as \eqref{eq5.6} with same upper-bounded constant $\mathcal{C}_3$.\\ \textbf{Case 2:} Now let's consider the case when  $\rho_3 (x)\equiv1 $ and $0<\rho_1,\rho_2\in L_1(\mathbb{R}_+)$. By the combination of \eqref{eq3.1} and \eqref{eq3.3} (for $\mathcal{I}_4$) together with Fubini's theorem, we obtain $$\begin{aligned}
\big|\underset{F_s,F_c,K}{\ast}(\rho_1,\rho_2,1)(x)\big| &\leq \frac{1}{2\sqrt{2\pi}}\int_{\mathbb{R}_+^3} |\Phi(x,u,v,w)| \rho_1(u)\rho_2(v) du dv dw\\&=\frac{1}{2\sqrt{2\pi}}\left(\int_0^\infty|\Phi(x,u,v,w)|dw\right)\left(\int_0^\infty \rho_1(u)du\right)\left(\int_0^\infty \rho_2(v)dv\right)\\& \leq \sqrt{\frac{2}{\pi}}~\mathrm{sech}\, t~\|\rho_1\|_{L_1(\mathbb{R}_+)}\|\rho_2\|_{L_1(\mathbb{R}_+)}<\infty.
\end{aligned}$$ This implies that operator $\underset{F_s,F_c,K}{\ast}(\rho_1,\rho_2,1)$ is well-defined and is a bounded
function  $\forall x\in\mathbb{R}_+$. Therefore, we have  
$\big|\underset{F_s,F_c,K}{\ast}(\rho_1,\rho_2,1)^{1-\frac1p}(x)\big| \leq \left(\sqrt{\frac{2}{\pi}}~\mathrm{sech}\,t\right)^{1-\frac1p} \|\rho_1\|^{1-\frac1p}_{L_1(\mathbb{R}_+)}\|\rho_2\|^{1-\frac1p}_{L_1(\mathbb{R}_+)}.$ Combining with the inequality \eqref{eq5.1}, we obtain the classical Saitoh's inequality as in \cite{Saitoh2000}
\begin{corollary}\label{cor5.2} Let  $0<\rho_1,\rho_2\in L_1(\mathbb{R}_+)$ and $\rho_3 (x) \equiv 1\ \forall x\in \R_+$. If $F_j\in L_p(\mathbb{R}_+,\rho_j),~j=1,2$ and  $F_3$ belongs to $L_p(\mathbb{R}_+)$, we have
	\begin{equation}\label{eq5.8}
	\big\|\underset{F_s,F_c,K}{\ast}(F_1\rho_1,F_2\rho_2,F_3)\big\|_{L_p(\mathbb{R}_+)}\leq \mathcal{C}_4\|F_1\|_{L_p(\mathbb{R}_+,\rho_1)}\|F_2\|_{L_p(\mathbb{R}_+,\rho_2)}\|F_3\|_{L_p(\mathbb{R}_+)}.
	\end{equation}
	where $\mathcal{C}_4=\sqrt{\frac{2}{\pi}}~(\mathrm{sech}\,t)^{1-\frac{1}{p}}\bigg[\sqrt{\frac{2}{\pi}}K_0(w)\bigg]^{\frac{1}{p}}\|\rho_1\|_{L_1 (\R_+)}^{1-\frac1p}\|\rho_2\|_{L_1 (\R_+)}^{1-\frac1p}$, with $p>1$.
\end{corollary}
\noindent Inequality \eqref{eq5.8} indicates that the mapping $\underset{F_s,F_c,K}{\ast}: L_p(\mathbb{R}_+,\rho_1)\times L_p(\mathbb{R}_+,\rho_2)\times L_p(\mathbb{R}_+) \to L_p(\mathbb{R}_+)$ defined by the impact $(F_1,F_2,F_3)\mapsto~\underset{F_s,F_c,K}{\ast}(F_1\rho_1,F_2\rho_2,F_3)$ is a bounded operator in $L_p(\mathbb{R}_+)$ for $p>1$. It is worth mentioning that all derived Saitoh-type inequalities in Subsection \ref{subsec32}, which include \eqref{eq5.1}, \eqref{eq5.6}, and \eqref{eq5.8} are always valid in $L_2(\mathbb{R}_+)$. Nevertheless, finding optimal sharp coefficients for these estimates remains challenging due to their dependence on weight functions.

\section{$L_1$-solution for classes of Toeplitz-Hankel type integral equation}\label{sec4}
\noindent This section will be devoted to classes of the Toeplitz-Hankel type integral equation related to \eqref{eq3.1}. Following \cite{Tsitsiklis1981Levy}, we consider the integral equation of the form \begin{equation}\label{eq7.1}
f(x)+\int\limits_{0}^{T} \big[k_1 (x+y)+k_2 (x-y) \big]f(y)dy=\vp (x), \ 0\leq x,y \leq T,
\end{equation}where $\vp$ is a given function, $f$ an  unknown function and $K(x,y):= k_1 (x+y)+k_2 (x-y)$ is the kernel. Equation \eqref{eq7.1} with a Hankel $k_1( x + y )$ or Toeplitz $k_2 ( x - y )$ kernel has attracted the attention of many authors due to its practical applications in diverse fields such as scattering theory, fluid dynamics, linear filtering theory, inverse scattering problems in quantum mechanics, and issues in radiative wave transmission, along with further applications in medicine and biology 
\cite{Agranovich1963Marchenko,Chadan1977Sabatier}.  
Eq. \eqref{eq7.1} has been carefully studied when $K(x,y)$ is a Toeplitz $k_2 (x - y)$ or Hankel kernel $k_1 (x + y)$. One notable case is when $K (x,y) = k_2 (x - y)+ k_2 ( x + y )$, i.e., the Toeplitz and Hankel kernels generated by the same function $k_2$ have been investigated in \cite{kagiwada1974integral}, specifically by establishing \eqref{eq7.1} in the finite interval $I = [0,T]$ with the kernel $K(x, \tau)=\int_{0}^{1}\big[e^{-\frac{|x-\tau|}{\theta}}+e^{-\frac{x+\tau}{\theta}} r(\theta)\big] w(\theta)\ d \theta.$ Tsitsiklis and Levy demonstrated the solvability of \eqref{eq7.1} in \cite{Tsitsiklis1981Levy} for general Toeplitz plus Hankel kernels $k_1 (x+y)+k_2 (x-y)$. This approach leads equation \eqref{eq7.1} to be a generalization
of Levinson's equation, as considered by Chanda-Sabatier-Newton in \cite{Chadan1977Sabatier} for the case where the kernel is a Toeplitz function, based on the Gelfand-Levitan's method \cite{Gelfand1951Levitan}. Similarly, and the approach by Marchenko \cite{Agranovich1963Marchenko} applies to the case where the kernel is a Hankel function. A special case of \eqref{eq7.1} encompasses a more extensive result of Tsitsiklis-Levy involving the KL-method with two Toeplitz plus Hankel-type kernels defined in the infinite domain $\mathbb{R}_+$ as follows: $$\lambda \int_{\mathbb{R}_{+}} K_1(x, \tau) f(\tau) d \tau+\mu \int_{\mathbb{R}_{+}} K_2(x, \tau) f(\tau) d \tau=g(x),$$ where the kernels take the form $$\begin{aligned}
	K_1(x, \tau)&=\frac{1}{\pi^2} \int_{\mathbb{R}_{+}}\big[\sinh (\tau-\theta) e^{-x \cosh (\tau-\theta)}+\sinh (\tau+\theta) e^{-x \cosh (\tau+\theta)}\big] \varphi(\theta) d \theta \quad \text{and}\\
K_2(x, \tau)&=\frac{1}{2 \pi x} \int_{\mathbb{R}_{+}}\big[e^{-x \cosh (\tau-\theta)}+e^{-x \cosh (\tau+\theta)}\big] h(\theta) d \theta,	
\end{aligned}$$  with complex constants $\lambda, \mu$ predetermined, as proven in \cite{tuan2018hoanghong}.
Being different from other approaches, our idea is to reduce the original integral equation to a linear equation by using the structure of polyconvolution \eqref{eq3.1}. Indeed,
we study the equation \eqref{eq7.1} by choosing the kernel pair $k_1$, $k_2$ as follows: $k_1(x+y)=-\frac{1}{\sqrt{2\pi}}\varphi(x+y)$; $k_2(x-y)=\frac{1}{\sqrt{2\pi}}\varphi(|x-y|)$ and $g(x)=\underset{F_s,F_c,K}{\ast}(\xi,\varphi,h)(x)$,  considering the infinite range $(0,T)\equiv (0,\infty)$, then the Toeplitz-Hankel type equation \eqref{eq7.1} can be rewritten in the convolution form 
\begin{equation}\label{eq7.2}
f(x)+(f\underset{1}{\ast} \varphi)(x)=\underset{F_s,F_c,K}{\ast}(\xi,\varphi,h)(x),\quad x>0.
\end{equation}
\begin{theorem}\label{thm7.1}
Let $\xi, \varphi,h$ be the given functions such that	$\xi, \varphi\in L_1(\mathbb{R}_+)$ and $h\in L_1^{0,\beta}(\mathbb{R}_+)$ with $\beta \in (0,1]$. Then, for the solvability of equation \eqref{eq7.2}  in $L_1 (\R_+)$, those are the sufficient condition including $1+(F_c\varphi)(y)\ne 0$. Moreover, its solution $f(x)$ is unique and represented by  $
f(x)=\underset{F_s,F_c,K}{\ast}(\xi,\ell,h)(x)=\big[(\xi\underset{1}{\ast} \ell)\underset{2}{\ast}h\big]
$ almost everywhere in $\R_+$, where $\ell\in L_1(\mathbb{R}_+)$ is defined via
$
(F_c\ell)(y)=\frac{(F_c\varphi)(y)}{1+(F_c\varphi)(y)}.
$ Here $(\underset{1}{\ast})$ is Sneddon's convolution \eqref{eq2.9} and $(\underset{2}{\ast})$ is Yakubovich-Britvina's convolution \eqref{eq2.12}.   Finally, the following $L_1$-norm estimation holds 
$\|f\|_{L_1(\mathbb{R}_+)}\leq \sqrt{\frac{2}{\pi}}\|\xi\|_{L_1(\mathbb{R}_+)}\|\ell\|_{L_1(\mathbb{R}_+)}\|h\|_{L_1^{0,\beta}(\mathbb{R}_+)}.$
\end{theorem}
\begin{proof}
Applying Fourier sine $(F_s)$ transform on both sides of \eqref{eq7.2}, we have	$$(F_sf)(y)+F_s(f\underset{1}{\ast}\varphi)(y)=F_s(\underset{F_s,F_c,K}{\ast}(\xi,\varphi,h))(y).$$ Based on the factorization property \eqref{eq3.6}, we obtain $(F_sf)(y)+(F_sf)(y)(F_c\varphi)(y)=(F_s\xi)(y)(F_c\varphi)(y)K[h](y)$, this equivalent to $(F_sf)(y)\big[1+(F_c\varphi)(y)\big]=(F_s\xi)(y)(F_c\varphi)(y)K[h](y)$. Under the condition $1+(F_c\varphi)(y)\ne 0, \forall y>0,$ it implies that $(F_sf)(y)=\frac{(F_c\varphi)(y)}{1+(F_c\varphi)(y)}(F_s\xi)(y)K[h](y)$ is valid. Applying the
Wiener–Lévy's theorem \cite{NaimarkMA1972} for the Fourier cosine $(F_c)$ transform,  there exists a function $\ell \in L_1(\mathbb{R}_+)$ such that $
(F_c\ell)(y)=\frac{(F_c\varphi)(y)}{1+(F_c\varphi)(y)},
\forall y>0$. This  yields $
(F_sf)(y)=(F_s\xi)(y)(F_c\ell)(y) K[h](y)=F_s\big(\underset{F_s,F_c,K}{\ast}(\xi,\ell,h)\big)(y).
$ Therefore, $f(x)=\underset{F_s,F_c,K}{\ast}(\xi,\ell,h)(x)$ almost everywhere $\forall x\in \mathbb{R}_+$. It is easy to point out $f(x) \in L_1 (\R_+)$. Indeed, 
	since $\xi\in L_1(\mathbb{R}_+)$, $h\in L_1^{0,\beta}(\mathbb{R}_+)$ and  $\ell\in L_1 (\R_+)$ due to Wiener–Lévy’s theorem, by Theorem \ref{thm3.1},
	we infer that $\underset{F_s,F_c,K}{\ast}(\xi,\ell,h)$ is  well-defined as a continuous function and belongs to $L_1 (\R_+)$. Moreover, due to Remark \ref{rmk1}, we obtain $f(x)=\underset{F_s,F_c,K}{\ast}(\xi,\ell,h)(x)=\big[(\xi \underset{1}{\ast} \ell)\underset{2}{\ast} h\big](x)$. The combination of these is sufficient to show there exists the only existence of $f$ belonging to $L_1(\mathbb{R}_+)$ because convolutions $(\underset{1}{\ast})$ and $(\underset{2}{\ast})$ are uniquely determined. Finally, we explicitly obtain this theorem's final assertion through estimate \eqref{eq3.7}.
\end{proof}

\noindent	It is important to emphasize that the condition stated in Theorem~\ref{thm7.1} is valid, ensuring that the existence of the function \( \ell \) through the Wiener--Lévy theorem is well-defined under the given assumptions. To elucidate this further, consider the following illustrative example. Let \( \varphi(x) = \sqrt{\frac{\pi}{2}} e^{-x} \). It is straightforward to verify that \( \varphi \in L_1(\mathbb{R}_+) \), and the  Fourier-cosine transform of \( \varphi \) is given by
$
	F_c\left(\sqrt{\frac{\pi}{2}} e^{-x}\right)(y) = \frac{1}{1 + y^2}.
$
	Consequently,
$
	1 + (F_c \varphi)(y) = 1 + \frac{1}{1 + y^2} \neq 0.
	$
	By invoking the Wiener--Lévy theorem, it follows that
	$
	(F_c \ell)(y) = \frac{(F_c \varphi)(y)}{1 + (F_c \varphi)(y)} = \frac{1}{2 + y^2}.
	$
	Referring to formula 1.4.1, p. 14 in \cite{bateman1954}, we conclude that 
	$
	\ell(x) = \frac{\sqrt{\pi}}{2} e^{-\sqrt{2}x} \in L_1(\mathbb{R}_+).
	$ Finally, applying inequalities \eqref{eq4.10} and \eqref{eq4.13} from Section~\ref{sec3}, we obtain an alternative estimate for the solution of equation~\eqref{eq7.2} as follows:
\begin{remark}\label{rem3}
\textup{A)	If $p,q,r,s$ be real numbers  in the open interval $(1,\infty)$ such that $\frac1p+\frac1q+\frac1r=2+\frac1s$. Assume that $ f\in L_s(\mathbb{R}_+)$, $\xi \in L_p(\mathbb{R}_+)$, $\ell \in L_q(\mathbb{R}_+)$ and $h\in L_r^{0,\beta}(\mathbb{R}_+)$, $\beta \in (0,1]$, then we obtain the solution’s estimation for the problem \eqref{eq7.2} as follows:
$
\|f\|_{L_s(\mathbb{R}_+)}\leq \mathcal{C}_1 \|\xi\|_{L_p(\mathbb{R}_+)}\|\ell\|_{L_q(\mathbb{R}_+)}\|h\|_{L_r^{0,\beta}(\mathbb{R}_+)}.
$	
Besides, for the case $s=\infty$ the solution of \eqref{eq7.2} is bounded $\forall x\in \mathbb{R}_+$ and we get the estimation}
	$
	\|f\|_{L_{\infty}(\mathbb{R}_+)}\leq \mathcal{C}_1 \|\xi\|_{L_p(\mathbb{R}_+)} \|\ell\|_{L_q(\mathbb{R}_+)} \|h\|_{L_r^{0,\beta}(\mathbb{R}_+)}, 
	$ \textup{where} $\mathcal{C}_1=\sqrt{\frac{2}{\pi}}\left(\mathrm{sech}\, t\right)^{1-\frac1r}$.\\
\textup{B) If the weighted function $\rho_1=1$,  without loss of generality, it can be assumed that the given functions $\ell=\ell_1 \rho_2$ and $h=h_1
	\rho_3$ such that $0<\rho_2,\rho_3\in L_1(\mathbb{R}_+)$, where $\ell_1\in L_1(\mathbb{R}_+,\rho_2)\cap L_p(\mathbb{R}_+,\rho_2)$, $h_1\in L_1(\mathbb{R}_+,\rho_3)\cap L_p(\mathbb{R}_+,\rho_3)$ and $f\in L_p(\mathbb{R}_+)$ with $p>1$. Using \eqref{eq5.6} provides us with an estimate of the boundedness in weighted
	$L_p$ spaces for the solution of equation \eqref{eq7.2} as follows} $
\|f\|_{L_p(\mathbb{R}_+)}\leq \mathcal{C}_3\|\xi\|_{L_p(\mathbb{R}_+)}\|\ell_1\|_{L_p(\mathbb{R}_+,\rho_2)}\|h_1\|_{L_p(\mathbb{R}_+,\rho_3)}$,  \textup{where the coefficient}
$
\mathcal{C}_3=\sqrt{\frac{2}{\pi}}K_0(w)\|\rho_2\|_{L_1(\mathbb{R}_+)}^{1-\frac1p}\|\rho_3\|_{L_1(\mathbb{R}_+)}^{1-\frac1p}.$
\end{remark}

\noindent \textbf{Disclosure statement}\\
\noindent No potential conflict of interest was reported by the author.\\
\noindent \textbf{Ethics declarations}\\
\noindent We confirm that all the research meets ethical guidelines and adheres to the legal requirements of the study country.
\noindent \textbf{Funding}\\
\noindent This research received no specific grant from any funding agency.\\
\noindent \textbf{ORCID}\\
\noindent \textsc{Trinh Tuan} {\color{blue} \url{https://orcid.org/0000-0002-0376-0238}}

\bibliographystyle{plain}
\bibliography{ref_ITSF_2025}

\end{document}